\documentclass[11pt]{amsart}
\usepackage{amsmath}
\usepackage{amssymb,latexsym}
\usepackage{amscd}
\usepackage{xspace}
\usepackage{verbatim, esint}
\begin{document}
\numberwithin{equation}{section}
\title[]{Semilinear elliptic equations admitting similarity transformations}
\author{Mousomi Bhakta, Moshe Marcus}
\address{Department of Mathematics, Technion\\
 Haifa, ISRAEL}
\date{\today}

\newcommand{\txt}[1]{\;\text{ #1 }\;}
\newcommand{\tbf}{\textbf}
\newcommand{\tit}{\textit}
\newcommand{\tsc}{\textsc}
\newcommand{\trm}{\textrm}
\newcommand{\mbf}{\mathbf}
\newcommand{\mrm}{\mathrm}
\newcommand{\bsym}{\boldsymbol}
\newcommand{\scs}{\scriptstyle}
\newcommand{\sss}{\scriptscriptstyle}
\newcommand{\txts}{\textstyle}
\newcommand{\dsps}{\displaystyle}
\newcommand{\fnz}{\footnotesize}
\newcommand{\scz}{\scriptsize}
\newcommand{\be}{\begin{equation}}
\newcommand{\bel}[1]{\begin{equation}\label{#1}}
\newcommand{\ee}{\end{equation}}
\newtheorem{subn}{\name}
\renewcommand{\thesubn}{}
\newcommand{\bsn}[1]{\def\name{#1$\!\!$}\begin{subn}}
\newcommand{\esn}{\end{subn}}
\newtheorem{sub}{\name}[section]
\newcommand{\dn}[1]{\def\name{#1}}   
\newcommand{\bs}{\begin{sub}}
\newcommand{\es}{\end{sub}}
\newcommand{\bsl}[1]{\begin{sub}\label{#1}}
\newcommand{\bth}[1]{\def\name{Theorem}\begin{sub}\label{t:#1}}
\newcommand{\blemma}[1]{\def\name{Lemma}\begin{sub}\label{l:#1}}
\newcommand{\bcor}[1]{\def\name{Corollary}\begin{sub}\label{c:#1}}
\newcommand{\bdef}[1]{\def\name{Definition}\begin{sub}\label{d:#1}}
\newcommand{\bprop}[1]{\def\name{Proposition}\begin{sub}\label{p:#1}}
\newcommand{\bnote}[1]{\def\name{\mdseries\scshape Notation}\begin{sub}\label{n:#1}}
\newcommand{\bproof}{\begin{proof}}
\newcommand{\eproof}{\end{proof}}
\newcommand{\bcom}{}
\newcommand{\req}{\eqref}
\newcommand{\rth}[1]{Theorem~\ref{t:#1}}
\newcommand{\rlemma}[1]{Lemma~\ref{l:#1}}
\newcommand{\rcor}[1]{Corollary~\ref{c:#1}}
\newcommand{\rdef}[1]{Definition~\ref{d:#1}}
\newcommand{\rprop}[1]{Proposition~\ref{p:#1}}
\newcommand{\rnote}[1]{Notation~\ref{n:#1}}
\newcommand{\BA}{\begin{array}}
\newcommand{\EA}{\end{array}}
\newcommand{\BAN}{\renewcommand{\arraystretch}{1.2}
\setlength{\arraycolsep}{2pt}\begin{array}}
\newcommand{\BAV}[2]{\renewcommand{\arraystretch}{#1}
\setlength{\arraycolsep}{#2}\begin{array}}
\newcommand{\BSA}{\begin{subarray}}
\newcommand{\ESA}{\end{subarray}}
\newcommand{\BAL}{\begin{aligned}}
\newcommand{\EAL}{\end{aligned}}
\newcommand{\BALG}{\begin{alignat}}
\newcommand{\EALG}{\end{alignat}}
\newcommand{\BALGN}{\begin{alignat*}}
\newcommand{\EALGN}{\end{alignat*}}
\newcommand{\note}[1]{\noindent\textit{#1.}\hspace{2mm}}
\newcommand{\Remark}{\textit{Remark}\hskip 2mm}
\newcommand{\forevery}{\quad \forall}
\newcommand{\1}{\\[1mm]}
\newcommand{\2}{\\[2mm]}
\newcommand{\3}{\\[3mm]}
\newcommand{\set}[1]{\{#1\}}
\def\capa{\mathit{cap}}
\newcommand{\st}[1]{{\rm (#1)}}
\newcommand{\lra}{\longrightarrow}
\newcommand{\lla}{\longleftarrow}
\newcommand{\llra}{\longleftrightarrow}
\newcommand{\Lra}{\Longrightarrow}
\newcommand{\Lla}{\Longleftarrow}
\newcommand{\Llra}{\Longleftrightarrow}
\newcommand{\warrow}{\rightharpoonup}
\def\dar{\downarrow}
\def\uar{\uparrow}
\newcommand{\paran}[1]{\left (#1 \right )}
\newcommand{\sqrbr}[1]{\left [#1 \right ]}
\newcommand{\curlybr}[1]{\left \{#1 \right \}}
\newcommand{\absol}[1]{\left |#1\right |}
\newcommand{\norm}[1]{\left \|#1\right \|}
\newcommand{\angbr}[1]{\left< #1\right>}
\newcommand{\paranb}[1]{\big (#1 \big )}
\newcommand{\sqrbrb}[1]{\big [#1 \big ]}
\newcommand{\curlybrb}[1]{\big \{#1 \big \}}
\newcommand{\absb}[1]{\big |#1\big |}
\newcommand{\normb}[1]{\big \|#1\big \|}
\newcommand{\angbrb}[1]{\big\langle #1 \big \rangle}
\newcommand{\thkl}{\rule[-.5mm]{.3mm}{3mm}}
\newcommand{\thknorm}[1]{\thkl #1 \thkl\,}
\newcommand{\trinorm}[1]{|\!|\!| #1 |\!|\!|\,}
\newcommand{\vstrut}[1]{\rule{0mm}{#1}}
\newcommand{\rec}[1]{\frac{1}{#1}}
\newcommand{\opname}[1]{\mathrm{#1}\,}
\newcommand{\supp}{\opname{supp}}
\newcommand{\dist}{\opname{dist}}
\newcommand{\sign}{\opname{sign}}
\newcommand{\diam}{\opname{diam}}
\newcommand{\q}{\quad}
\newcommand{\qq}{\qquad}
\newcommand{\hsp}[1]{\hspace{#1mm}}
\newcommand{\vsp}[1]{\vspace{#1mm}}
\newcommand{\prt}{\partial}
\newcommand{\sms}{\setminus}
\newcommand{\ems}{\emptyset}
\newcommand{\ti}{\times}
\newcommand{\pr}{^\prime}
\newcommand{\ppr}{^{\prime\prime}}
\newcommand{\tl}{\tilde}
\newcommand{\wtl}{\widetilde}
\newcommand{\sbs}{\subset}
\newcommand{\sbeq}{\subseteq}
\newcommand{\nind}{\noindent}
\newcommand{\ovl}{\overline}
\newcommand{\unl}{\underline}
\newcommand{\nin}{\not\in}
\newcommand{\pfrac}[2]{\genfrac{(}{)}{}{}{#1}{#2}}
\newcommand{\tin}{\to\infty}
\newcommand{\ind}[1]{_{_{#1}}}
\newcommand{\chr}[1]{\mathbf{1}\ind{#1}}
\newcommand{\rest}[1]{\big |\ind{#1}}
\newcommand{\Sol}[2]{\mathrm{Sol}\ind{#2}{#1}}
\newcommand{\wkc}{weak convergence\xspace}
\newcommand{\wrto}{with respect to\xspace}
\newcommand{\cons}{consequence\xspace}
\newcommand{\consy}{consequently\xspace}
\newcommand{\Consy}{Consequently\xspace}
\newcommand{\Essy}{Essentially\xspace}
\newcommand{\essy}{essentially\xspace}
\newcommand{\mnz}{minimizer\xspace}
\newcommand{\sth}{such that\xspace}
\newcommand{\ngh}{neighborhood\xspace}
\newcommand{\nghs}{neighborhoods\xspace}
\newcommand{\seq}{sequence\xspace}
\newcommand{\seqs}{sequences\xspace}
\newcommand{\sseq}{subsequence\xspace}
\newcommand{\ifif}{if and only if\xspace}
\newcommand{\suff}{sufficiently\xspace}
\newcommand{\abc}{absolutely continuous\xspace}
\newcommand{\sol}{solution\xspace}
\newcommand{\subss}{sub-solutions\xspace}
\newcommand{\subs}{sub-solution\xspace}
\newcommand{\supers}{super-solution\xspace}
\newcommand{\superss}{super-solutions  \xspace}
\newcommand{\Wlg}{Without loss of generality\xspace}
\newcommand{\wlg}{without loss of generality\xspace}
\newcommand{\locun}{locally uniformly\xspace}
\newcommand{\bvp}{boundary value problem\xspace}
\newcommand{\bvps}{boundary value problems\xspace}
\newcommand{\bdw}{\partial\Gw}
\newcommand{\Capq}{C_{2/q,q'}}
\newcommand{\Wq}{W^{2/q,q'}}
\newcommand{\Wqdual}{W^{-2/q,q}}
\newcommand{\Wqdb}{W^{-2/q,q}_+(\bdw)}
\newcommand{\sbsq}{\overset{q}{\sbs}}
\newcommand{\smq}{\overset{q}{\sim}}
\newcommand{\app}[1]{\underset{#1}{\approx}}
\newcommand{\suppq}{\mathrm{supp}^q_{\bdw}\,}
\newcommand{\convq}{\overset{q}{\to}}
\newcommand{\barq}[1]{\bar{#1}^{^q}}
\newcommand{\prtq}{\partial_q}
\newcommand{\tr}{\mathrm{tr}\,}
\newcommand{\Tr}{\mathrm{Tr}\,}
\newcommand{\trin}[1]{\mathrm{tr}\ind{#1}}
\newcommand{\Ltrin}[1]{\text{\rm L-tr}\ind{#1}}
\newcommand{\Mtrin}[1]{\text{\rm M-tr}\ind{#1}}
\newcommand{\qcl}{$q$-closed\xspace}
\newcommand{\qop}{$q$-open\xspace}
\newcommand{\gsmod}{$\gs$-moderate\xspace}
\newcommand{\gsreg}{$\gs$-regular\xspace}
\newcommand{\qreg}{quasi regular\xspace}
\newcommand{\qeq}{$q$-equivalent\xspace}
\newcommand{\ppf}{\underset{f}{\prec\prec}}
\newcommand{\ofrown}{\overset{\frown}}
\newcommand{\modcon}{\underset{mod}{\lra}}
\newcommand{\ugb}[1]{u\chr{\Gs_\gb(#1)}}
\newcommand{\mcon}{$q$-moderately convergent\xspace}
\newcommand{\mdiv}{$q$-moderately divergent\xspace}
\def\qsupp{q\text{-supp}\,}
\def\RN{\BBR^N}
\def\loc{\ind{\rm loc}}
\def\bcom{}
\def\ga{\alpha}     \def\gb{\beta}       \def\gg{\gamma}
\def\gc{\chi}       \def\gd{\delta}      \def\ge{\epsilon}
\def\gth{\theta}                         \def\vge{\varepsilon}
\def\gf{\phi}       \def\vgf{\varphi}    \def\gh{\eta}
\def\gi{\iota}      \def\gk{\kappa}      \def\gl{\lambda}
\def\gm{\mu}        \def\gn{\nu}         \def\gp{\pi}
\def\vgp{\varpi}    \def\gr{\rho}        \def\vgr{\varrho}
\def\gs{\sigma}     \def\vgs{\varsigma}  \def\gt{\tau}
\def\gu{\upsilon}   \def\gv{\vartheta}   \def\gw{\omega}
\def\gx{\xi}        \def\gy{\psi}        \def\gz{\zeta}
\def\Gg{\Gamma}     \def\Gd{\Delta}      \def\Gf{\Phi}
\def\Gth{\Theta}
\def\Gl{\Lambda}    \def\Gs{\Sigma}      \def\Gp{\Pi}
\def\Gw{\Omega}     \def\Gx{\Xi}         \def\Gy{\Psi}

\def\CS{{\mathcal S}}   \def\CM{{\mathcal M}}   \def\CN{{\mathcal N}}
\def\CR{{\mathcal R}}   \def\CO{{\mathcal O}}   \def\CP{{\mathcal P}}
\def\CA{{\mathcal A}}   \def\CB{{\mathcal B}}   \def\CC{{\mathcal C}}
\def\CD{{\mathcal D}}   \def\CE{{\mathcal E}}   \def\CF{{\mathcal F}}
\def\CG{{\mathcal G}}   \def\CH{{\mathcal H}}   \def\CI{{\mathcal I}}
\def\CJ{{\mathcal J}}   \def\CK{{\mathcal K}}   \def\CL{{\mathcal L}}
\def\CT{{\mathcal T}}   \def\CU{{\mathcal U}}   \def\CV{{\mathcal V}}
\def\CZ{{\mathcal Z}}   \def\CX{{\mathcal X}}   \def\CY{{\mathcal Y}}
\def\CW{{\mathcal W}}
\def\BBA {\mathbb A}   \def\BBb {\mathbb B}    \def\BBC {\mathbb C}
\def\BBD {\mathbb D}   \def\BBE {\mathbb E}    \def\BBF {\mathbb F}
\def\BBG {\mathbb G}   \def\BBH {\mathbb H}    \def\BBI {\mathbb I}
\def\BBJ {\mathbb J}   \def\BBK {\mathbb K}    \def\BBL {\mathbb L}
\def\BBM {\mathbb M}   \def\BBN {\mathbb N}    \def\BBO {\mathbb O}
\def\BBP {\mathbb P}   \def\BBR {\mathbb R}    \def\BBS {\mathbb S}
\def\BBT {\mathbb T}   \def\BBU {\mathbb U}    \def\BBV {\mathbb V}
\def\BBW {\mathbb W}   \def\BBX {\mathbb X}    \def\BBY {\mathbb Y}
\def\BBZ {\mathbb Z}

\def\GTA {\mathfrak A}   \def\GTB {\mathfrak B}    \def\GTC {\mathfrak C}
\def\GTD {\mathfrak D}   \def\GTE {\mathfrak E}    \def\GTF {\mathfrak F}
\def\GTG {\mathfrak G}   \def\GTH {\mathfrak H}    \def\GTI {\mathfrak I}
\def\GTJ {\mathfrak J}   \def\GTK {\mathfrak K}    \def\GTL {\mathfrak L}
\def\GTM {\mathfrak M}   \def\GTN {\mathfrak N}    \def\GTO {\mathfrak O}
\def\GTP {\mathfrak P}   \def\GTR {\mathfrak R}    \def\GTS {\mathfrak S}
\def\GTT {\mathfrak T}   \def\GTU {\mathfrak U}    \def\GTV {\mathfrak V}
\def\GTW {\mathfrak W}   \def\GTX {\mathfrak X}    \def\GTY {\mathfrak Y}
\def\GTZ {\mathfrak Z}   \def\GTQ {\mathfrak Q}
\font\Sym= msam10
\def\SYM#1{\hbox{\Sym #1}}
\def\rqq{\req{q-eq}\xspace}
\def\bmn{\mathbf{n}}
\def\bma{\mathbf{a}}
\newcommand{\prtn}{\prt_{\bmn}}
\def\txin{\txt{in}}
\def\txon{\txt{on}}

\maketitle
\begin{center}{\small Department of Mathematics, Technion,
Haifa 32000, ISRAEL, \\
email:  \textit{mousomi@techunix.technion.ac.il, \ marcusm@math.technion.ac.il}}
\end{center}
\begin{abstract}
In this paper we study the equation $-\Gd u+\rho^{-(\ga+2)}h(\rho^{\ga}u)=0$ in a smooth bounded domain $\Gw$ where $\rho(x)=\text{dist}(x,\bdw), \q\ga>0$ and $h$ is a nondecreasing function which satisfies Keller-Osserman condition.  We introduce a  condition on $h$  which implies that the equation is subcritical, i.e. the corresponding boundary value problem is well posed with respect to data given by  finite measures. Under additional assumptions on $h$ we show that this condition is necessary as well as sufficient. We also discuss b.v. problems with data given by positive unbounded measures. Our results extend results of \cite{MV1} treating equations of the form $-\Delta u+\rho^\beta u^q=0$ with $q>1$, $\beta>-2$.
\vskip 2mm

\noindent keywords:\hskip 2mm\textit{semilinear elliptic equations, similarity transformation, very singular solution, isolated singularity, generalized trace, generalized boundary value problem, measure data.}
\vskip 2mm

\end{abstract}

\section{\bf Introduction}

In this article we consider  equations of the type
\begin{equation}\label{1}
-\Gd u+H(\rho, u)=0  \quad \text{in} \quad \Gw,
\end{equation}
where  $\Gw$ is a $C^2$ domain in $\BBR^N$, $N>1$,
\begin{equation}\label{def-H}
H(\rho, u)=\rho^{-2-\ga}h(\rho^\ga u),\q\rho(x)=\dist(x,\bdw)
\end{equation}
and  $\ga>0$. With respect to the nonlinearity $h$ we assume:
\begin{equation}\label{2}
\BAL
&(i)&& h\in C^1(\BBR),\ h(0)=0, \ h(t)>0 \ \text{for} \ t>0,\\
&(ii)&& h \ \text{is a convex and odd function},\\
&(iii)&& \int_{\Gw\cap B_R(0)}h(c\rho^\ga)\rho^{-(1+\ga)}dx<\infty \forevery c>0,\;R>0.
\EAL
\end{equation}
Note that (i) and (ii) imply that $h$ is nondecreasing.

More general equations such as
\begin{equation}\label{1'}
-\Gd u+g\circ u=0  \quad \text{in} \quad \Gw,
\end{equation}
(where $g\circ u (x)=g(x, u(x))$) and various special cases have been studied intensively over the last twenty years (see \cite{MV}, \cite{MV1} and the references therein). In the case of equations with absorption, a basic set of assumptions on $g$ is
\begin{equation}\label{2'}
\BAL
&(i) && g(x, \cdot)\in C^1(\BBR),\ g(x, 0)=0, \ g(x, t)>0 \ \text{for} \ t>0,\\
&(ii) && g(x,\cdot) \ \text{is nondecreasing and odd},\\
&(iii) && g(\cdot, t)\in L^{1}(\Gw\cap B_R(0),\rho) \ \forevery \ t\in\BBR \ \text{and}\ \forevery\ R>0,
\EAL
\end{equation}
 The assumption that $g(x,\cdot)$ is odd  is often omitted. However it simplifies the presentation without affecting the treatment in an essential way. We note that conditions \req{2} imply that $g=H$ satisfies \req{2'}.

The family of functions satisfying  \eqref{2'} will be denoted by
$\mathcal{G}_0=\mathcal{G}_0(\Gw)$.  Equations of the form \req{1} (with $\gr(x)$ replaced by $|x|$) have been introduced by Bandle and Marcus \cite{BM1}. These equations admit a similarity transformation. A special case of \req{1}, namely,
 $$H(\gr,t)=\gr^\gg |t|^q\sign t, \q \gg>-2$$
 has been extensively studied. The case $\gg=0$ received much attention;
    the problem with arbitrary $\gg>-2$ was studied in \cite{MV1}. For the case $\gg>0$ see \cite{DG} and references therein. Equation \eqref{def-H} includes the special case $H(\rho, t)=t^p+\rho^{\gg}t^q$ where $p,q>1, \gg>-2$.

\medskip

\noindent\textit{Notation.}\hskip 2mm  We denote by $\CM(\bdw)$ the family of set functions $\nu$ on $\CB(\bdw)$, such that, for every compact set $E\sbs \bdw$, $\nu\chr{E}$ is a finite measure. Thus, if $\bdw$ is compact then $\CM(\bdw)$ is the set of finite Borel measures on $\bdw$.

Put
$$\tilde L^1(\Gw)=\{u\in L^1_{loc}(\Gw): u\in L^1(\Gw\cap B_R(0)) \ \forall R>0\}$$
and
$$C^2_0(\bar\Gw):=\{\phi\in C^2(\bar\Gw) : \phi=0 \ \text{on} \ \bdw\}.$$
Assume that $\Gw$ is a $C^2$ domain (not necessarily bounded) and let $\nu\in \CM(\bdw)$.

 A function $u$ ia a solution of \emph{equation} \eqref{1'} if $u$ and
$g\circ u$ are locally integrable in $\Gw$  and the equation is satisfied  in the distribution sense.

A function $u$ is a solution of the \emph{boundary value problem}
\begin{equation}\label{eq-g}
\BAL
-\Gd u+g\circ u &=0 \ \ &&\text{in} \ \Gw; \\
  u &=\nu \ \ &&\text{on} \ \bdw,
\EAL
\end{equation}
 where $\nu$ is a Radon measure on $\bdw$, if $u$ and  $(g\circ u)\rho$ belong to $\tilde L^1(\Gw)$ and
\begin{equation}\label{g-weak}
\int_{\Gw}(-u\Gd\phi+(g\circ u)\phi)dx=-\int_{\bdw}\frac{\partial\phi}{\partial{\bf n}}d\nu
\end{equation}
for every $\phi\in C^2_0(\bar\Gw)$ with bounded support.

The set of all measures $\nu\in\CM(\bdw)$ such that \eqref{eq-g} possesses a solution is denoted by $\CM^g(\bdw)$. It is well-known that when a solution exists it is unique, \cite{MV}.

\bdef{KO}
Let $f :\BBR\to\BBR$ be an odd function and satisfy the following assumptions:
$$f\in C^1(\BBR), \ \ f(0)=0, \ \ f(t)>0 \ \text{for} \ t>0 .$$
We say that $f$ satisfies Keller-Osserman condition (KO) if
\begin{equation}\label{3}
\int_{a}^{\infty}F(s)^{-\frac{1}{2}}ds<\infty\ \ \forall a>0,\ \ \text{where} \ F(s)=
\int_{0}^{s}f(t)dt.
\end{equation}

If $g(x,t)\in\mathcal{G}_0$ we say that it satisfies the
global KO condition if there exists $f$ as above that satisfies the KO condition and
$$f(|t|)\leq|g(x,t)| \quad x\in\Gw, \ t\in\BBR.$$

 We say that $g$ satisfies the local KO condition if, for every domain $\Gw^{'}\Subset\Gw$ there exists $f=f^{\Gw^{'}}\in C(\BBR)$ that satisfies the KO condition and
$$f(|t|)\leq|g(x,t)| \quad x\in\Gw^{'}, \ t\in\BBR.$$
\es

\nind{\bf \Remark \ref{d:KO}:} If $h$  satisfies the KO condition then the function $H$ given by \eqref{def-H} satisfies the {\it local} KO condition. Therefore, in this case, the family of solutions of \eqref{1} is uniformly bounded in compact subsets of $\Gw$.

\vspace{2mm}

\bdef{trace}Let $u\in W^{1,p}_{loc}(\Gw)$ for some $p>1$. We say that a measure $\nu\in\CM(\bdw)$ is the boundary trace of $u$ on $\bdw$ if, for every uniform $C^2$ exhaustion $\{{\Gw}_n\}$,
\begin{equation}\label{trace}
    \int_{{\bdw}_n} u\chr{\bdw_n}\phi dS\to\int_{\bdw}\phi d\nu
\end{equation}
for every compactly supported $\phi\in C(\bar\Gw)$.
(Here $u\chr{\bdw_n}$ denotes the Sobolev trace.)

Similarly we define the boundary trace on a relatively open set $A\sbs \bdw$. A measure $\nu\in\CM(A)$ is the boundary trace of $u$ on $A$ if  \req{trace} holds for every $\phi\in C(\bar\Gw)$ \sth $\supp\phi$ is a compact subset of $\Gw\cup A$. In the case of positive solutions we slightly extend this definition  to include positive Radon measures on $A$.
 \es

\blemma{reg-sol} Let $g\in \CG_0$ and let $u$ be a positive solution of \req{1'}. Let $O$ be an open set  in $\BBR^N$ \sth $O\cap\bdw\neq\ems$. If
\begin{equation}\label{reg-cond}
   \int_{\Gw\cap O}g(x,u(x))\rho(x)dx<\infty.
\end{equation}
Then $u$ has a boundary trace on $O\cap\bdw$.
\es

\nind\textit{Notation.}\hskip 2mm We denote by $\CB_{reg}(\bdw)$ the space of positive  regular Borel measure on $\bdw$.

 Note that a measure $\nu\in\CB_{reg}(\bdw)$ need not be a Radon measure; it may blow up on compact sets. However the outer regularity implies that, for each measure $\nu$ in this space there exists a closed set $S_\nu$, called \emph{ the blow up} set of $\nu$,
  \sth $\nu(A)=\infty$ for every non-empty Borel set $A\sbs S_\nu$ and the restriction of $\nu$ to $R_\nu:=\bdw\sms S_\nu$ is a Radon measure.

Next we extend the notion of boundary trace to positive solutions of \req{1'} that may not  have a finite boundary trace.
  \bdef{ext-tr} Let $u$ be a positive solution of \req{1'}. We say that $u$ has a (generalized) boundary trace $\nu\in \CB_{reg}(\bdw)$ if:

   (a)   $\nu\lfloor_{R_\nu}$ is the boundary trace of $u$ on the relatively open set $R_\nu$ (in the sense of \rdef{trace}) and

  (b) for every open set $O$ in $\BBR^N$\sth $O\cap S_\nu\neq\ems$ and every uniform $C^2$ exhaustion $\{{\Gw}_n\}$ (see \cite[Definition 1.3.1]{MV})
\begin{equation}\label{ext-tr}
    \int_{{\bdw}_n\cap O} u \,dS\to\infty.
\end{equation}
\es

\bdef{reg-pt}
Let $u\in C^{2}(\Gw)$ be a positive solution of \eqref{1'}. A point $z\in\bdw$ is called a  regular boundary point of $u$ if there exists an open neighborhood $U$ of $z$ such that
\begin{equation}\label{40}
\int_{U\cap\Gw}(g\circ u)\rho\,dx<\infty.
\end{equation}
The set of regular boundary points of $u$ will be denoted by $\mathcal{R}(u)$. Its complement $\CS(u):=\bdw\setminus\mathcal{R}(u)$ is the singular boundary set of $u$. A point $y\in \CS(u)$ is called a strongly singular point of $u$.
\es
Clearly $\CR(u)$ is relatively open; therefore $\CS(u)$ is closed.

Here we recall \cite{MV}\cite[Theorem 1.1]{MV1}.

\bth{existence4}
Let $\Gw$ be a domain whose boundary is a manifold of class $C^2$, not necessarily compact. Suppose that $g\in\mathcal{G}_0$, satisfies  the local KO condition. Also assume that \eqref{1'} possesses a {\it barrier} at every point of $\bdw$.  Then every positive solution $u$ of \eqref{1'} possesses a (generalized) boundary trace given by a positive measure $\nu\in \CB_{reg}(\bdw)$.

$\CS(u)$ coincides with the blow up set of $\nu$ so that $\nu$ is a Radon measure on $\CR(u)$.
\es

For the definition of { barrier} and the conditions of its existence  see Definition \ref{d:barrier} and Proposition \ref{p:barrier}.
\vskip 3mm

\bdef{subcriticality} A nonlinearity $g\in\mathcal{G}_0$ is called {\it subcritical} if  problem \eqref{eq-g} possesses a solution for every $\nu\in\CM(\bdw)$.
\es

In this article we focus on the boundary value problem
\begin{equation}\label{1''}
\BAL
-\Gd u+H(\rho, u)&=0  \quad &&\text{in} \quad \Gw;\\  u&=\nu \ \ &&\text{on} \ \ \bdw,
\EAL
\end{equation}
where $\nu$ is a regular Borel measure, $H$ is as in \eqref{def-H}  and $h$ is assumed to satisfy \eqref{2}.

\vspace{2mm}

Following is a description of our main results concerning this problem.

\medskip

\bth{continuity}
Let $\Gw$ be a bounded domain of class $C^2$ and  $H$ is defined as in \eqref{def-H}. Assume that $h$ satisfies \req{2}(i),(ii) and
\begin{equation}\label{subcon1}
H(\rho,cP(\cdot,y))\in L^1(\Gw;\rho) \forevery y\in \bdw, \  \forevery  c\in\BBR.
\end{equation}
where $P(x,y)$ is the Poisson kernel of $-\Gd$ in $\Gw$.

Then $H$ is subcritical and the following assertion holds:

Assume that  there exists $\ge>0$ such that
\begin{equation}\label{h-zero}
\limsup_{t\to 0}\frac{h(t)}{t^{1+\epsilon}}<\infty.
\end{equation}
Let $\{\nu_k\}$ be a \seq in $\CM(\bdw)$ converging weakly to a measure $\mu$. Let $v$ (resp. $v_k$) denote the solution of \req{1''} with $\nu=\mu$ (resp. $\nu=\nu_k$).
Then $v_k\to v$ in $L^1(\Gw)$ and $H(\rho, v_k)\to H(\rho, v)$ in $L^1(\Gw,\rho)$.
\es

\medskip

\noindent\textit{Remark \ref{t:continuity}} Let $y\in \bdw$ and put
$$C_R=\{x\in\Gw\cap B_R(y) : |x-y|<2\rho(x)\}.$$
Assume that $R$ is sufficiently small so that $\bar C_R\subset\Gw\cup\{y\}$. Then \req{subcon1} implies
$$\int_{C_R}\rho^{-(1+\ga)}h(c\rho^{\ga}P(x,y))dx<\infty$$
which in turn implies
\begin{equation}\label{sub}
\int_{0}^{1}t^{N-\ga-2}h(ct^{\ga-N+1})dt<\infty.
\end{equation}
Actually, in a bounded $C^2$ domain \req{sub} is equivalent to \req{subcon1} (see Section 4).
\eqref{sub} implies that, $\ga\neq N-1$ and
\begin{equation}\label{subcon2}\BAL
(a)&\quad\int_{a}^{\infty}h(t^{-1})dt<\infty \quad&&\text{if} \ \ga>N-1\\
(b)&\quad\int_{0}^{a}h(t^{-1})dt <\infty \quad&&\text{if} \ 0<\ga<N-1,
\EAL
\end{equation}
for every $a>0$.
 Consequently,  if $h$ satisfies \req{2}:
\begin{equation}\label{N-1}
 \textrm{Condition  (\ref{subcon1})  implies $\ga>N-1$.}
\end{equation}

\noindent Indeed, as  $h$ is convex on $[0,\infty)$, $h(0)=0$ and $h$ is nondecreasing, it follows that $s\mapsto h(s)/s$ is nondecreasing. Thus $h(s)\geq h(1)s$ for $s>1$ and, by assumption, $h(1)>0$. Therefore $\int_{0}^{1}h(t^{-1})dt=\infty$.  This rules out  \req{subcon2}(b) so that $\ga>N-1$.

 In Section 4 we show  that, under some additional assumptions on $h$, the condition $\ga>N-1$ is \textit{necessary and sufficient for $H$ to be subcritical}; in particular it is equivalent to the subcriticality condition \eqref{subcon1}.

\vspace{2mm}

If  $h(t)=t^q,\ q>1$ then $H=\rho^{\gb}t^q$ where $\gb=\ga(q-1)-2$. In this case, by \cite{MV1}, $H$ is subcritical if and only if $q<\frac{N+1+\gb}{N-1}$, i.e., $\ga>N-1$.
\vskip 2mm

\bdef{delta-2}We say that $h$ satisfies $\Gd_2$ condition if there exists a constant $c>0$ such that
$$h(a+b)\leq c\big(h(a)+h(b)\big) \ \forall a>0, \ b>0.$$
\es

\bth{removability} Let $\Gw$ be a bounded domain, $y\in\bdw$ and $H$ be defined as in \eqref{def-H}. Suppose that $h$ satisfies KO condition, $\Gd_2$ condition, \eqref{2} and $H$ satisfies the global barrier condition.
If  $\ga\leq N-1$ and $u$ is a nonnegative solution of \eqref{1}  vanishing on $\bdw\setminus\{y\}$, then $u=0$ in $\Gw$.
\es

\noindent\textit{Notation.}\hskip 2mm Assume that $H$ is subcritical. If $y\in \bdw$ denote by $u_{k,y}$ the unique solution of \eqref{1''} with $\nu=k\gd_{y}$.

\bth{lower-estimate}
Let $\Gw$ be a bounded $C^2$ domain. Suppose that $H$ satisfies \eqref{subcon1}  and the global barrier condition.  In addition assume  that  $h$ satisfies the KO condition, \eqref{2} and \eqref{h-zero}.

Under these conditions,
\begin{equation}\label{u_infty}
   u_{\infty,y}=\lim_{k\tin}u_{k,y}
\end{equation}
is a solution of \req{1} that vanishes on $\bdw\sms\{y\}$. Furthermore,
if $u$ is a positive solution of \eqref{1} and $y\in\mathcal{S}(u)$ then
\begin{equation}\label{u_infty<u}
     u_{\infty,y}\leq u\q\text{in }\Gw.
\end{equation}
\es
For the term 'global barrier condition', see Definition \ref{d:barrier}.

\bdef{iss} Let $y\in \bdw$. If $u$  is a positive solution of \req{1} \sth $u=0$ on $\bdw\sms \{y\}$ and $y\in \CS(u)$ we say that $u$ has a \emph{strong isolated singularity} at $y$.
\es

\rth{lower-estimate} implies that $u_{\infty,y}$ is the \emph{smallest} such solution.

For the next statement we need some additional notation. Let $\Gw$ be a $C^2$ domain and let $y\in \bdw$. We denote by $\mathbf{n}_y$ the unit normal to $\bdw$ at $y$ pointing outward. We denote by $\mathbf{R}_y$ the rotation that maps the vector $e_1=(1,0,\ldots,0)$ to $-\mathbf{n}_y$.

\bth{existence2}
Let $\Gw$ be a bounded $C^2$ domain ,  $y\in\bdw$. Suppose that $H$ and $h$ satisfy the assumptions of \rth{lower-estimate}. In addition assume that $h$ is strictly convex.

Then

(i)\q equation \eqref{1} possesses a unique solution with strong isolated singularity at $y$.  This solution, denoted by $U_y$, satisfies
\begin{equation}\label{60}
\lim_{r\to 0}r^{\ga} U_{y}(y+r\sigma)=w(\mathbf{R}_y(\sigma)),
\end{equation}
where $w$ is the solution to the problem
\begin{equation}\label{52}
\BAL
\Gd_{S^{N-1}}w+\gl w-(\sigma\cdot e_1)^{-(2+\ga)}h\big((\sigma\cdot e_1)^{\ga}w\big)&=0 \quad \text{in} \ \ S^{N-1}_{+},\\
w &= 0 \quad \text{on} \ \ \partial S^{N-1}_{+},
\EAL
\end{equation}
$\gl=\ga-(N-1)$, $e_1=(1,0,\ldots,0)$ and
$S^{N-1}_{+}=\{x\in S^{N-1}:x_1>0\}$.
The convergence  is uniform in compact subsets of $S^{N-1}_{+}$.

\vspace{2mm}

(ii)\q There exists a positive constant $C$, depending on $N,\ga, h$, $C^2$ characteristic of $\Gw$ but independent of $y$, such that
\begin{equation}\label{asymp-esti}
C^{-1}|x-y|^{-\ga-1}\rho(x)\leq U_y(x)\leq C^{1}|x-y|^{-\ga-1}\rho(x) \q\forall x\in \Gw.
\end{equation}
\es

\bth{existence3}
Let $\Gw$ be a bounded $C^2$ domain. Under the assumptions of \rth{existence2}, problem \eqref{1''} possesses a unique positive solution for every $\nu\in\mathcal{B}_{reg}(\bdw)$. Furthermore, if $F$ is a closed subset of $\bdw$ then the unique solution, $U_F$, with the boundary trace $\infty\CX_{F}$ satisfies,
\begin{equation}\label{asymp-estii}
 c^{-1} \text{dist}(x,F)^{-\ga-1}\rho(x)\leq U_F(x)\leq c \text{dist}(x,F)^{-\ga-1}\rho(x) \2 \q\forall x\in\Gw.
 \end{equation}
\es

Some examples of nonlinearities for which the above theorems hold:
\begin{equation}\label{examples1}\BAL
&(i)  &&H(x,t)=\rho^\beta |t|^q \sign t, \quad q>1,\;\beta>-2, \\
&(ii) &&H(x,t)= \rho^\beta  |t|^q \ln (1+\rho^\ga |t|)\sign t  \quad q>1,\;\beta=-2+\ga(q-1),\\
&(iii)&& H(x,t)=\rho^{-2}|t|\exp^{\rho^{\ga}|t|}\sign t.
\EAL\end{equation}

The first example was treated in \cite{MV1} and similar results have been obtained. However the estimates corresponding to \req{asymp-esti} and \req{asymp-estii} have been established only for $\beta\geq0$.

\vspace{2mm}

The paper is organized as follows: In section 2 we present some definitions and preliminary results.  In section 3 we derive global estimates for solutions of \eqref{1}. Section 4 is devoted to the proof of Theorem \ref{t:continuity} and Theorem \ref{t:removability}. We establish a lower estimate of the singular solution in Section 5. Section 6 contains the proof of existence and uniqueness of solutions with strong isolated singularity, first in half space and then in bounded domain. In section 7 we prove Theorem \ref{t:existence3}.

\section{\bf Definitions and preliminary results}

Let $\Gw$ be a $C^2$ domain in $\BBR^N$. We use our notation:
\begin{equation}\label{4}
\Gs_{\gb}=\{x\in\Gw:\rho(x)=\gb\}, \quad D_{\gb}=\{x\in\Gw:\rho(x)>\gb\}, \quad
\Gw_{\gb}=\Gw\setminus D_{\gb}.
\end{equation}
${\bf n}_{x}$ denotes the unit outward normal at $x\in\bdw$. It is known that there exists $\gb_0>0$ such that:
\begin{itemize}
\item[(a)] For every point $x\in\bar\Gw_{\gb_{0}}$, there exists a unique point $\gs(x)\in\bdw$ such that $|x-\gs(x)|=\rho(x)$, i.e.
$$x=\gs(x)-\rho(x){\bf n}_{\gs(x)}, \quad \lim_{x\to\gs(x)}\nabla\rho(x)=-{\bf n}_{\gs(x)}.$$
\item[(b)]The mapping $x\mapsto\rho(x)$ and $x\mapsto\gs(x)$ belong to $C^2(\bar\Gw_{\gb_{0}})$ and $C^1(\bar\Gw_{\gb_{0}})$ respectively.
\item[(c)] The mapping $x\mapsto (\rho(x),\gs(x))$ is a $C^2$ diffeomorphism from $\bar\Gw_{\gb_{0}}$ to $[0,\gb_0]\times\bdw$.
\end{itemize}
Thus $(\rho,\gs)$ may serve as a set of coordinate in one sided neighborhood of the boundary and are called \textit{flow coordinates} of $\Gw$.

\bdef{barrier} Let $z\in\bdw$. We say that there exists a g-barrier at $z$ if there exists $r(z)>0$ such that for every $0<r\leq r(z)$, \eqref{1'} possesses a positive solution
$w=w_{r,z}$ in $B_{r}(z)\cap\Gw$ such that
\begin{equation}\label{barrier}
\BAL
&(i)&& w\in C(B_{r}(z)\cap\bar\Gw),\quad w=0 \ \text{on} \ \bdw\cap B_{r}(z), \\
&(ii)&& {\lim_{x\to y}}w(x)=\infty \quad \forall y\in\partial B_{r}(z)\cap\Gw.
\EAL
\end{equation}
We say that $g$ satisfies the global barrier condition if:\\
(i) A g-barrier exists at every point of $\bdw$.\\
(ii) There exists a number $\bar r>0$ such that $r(z)\geq\bar r \quad \forall \ z\in\bdw$.\\
(iii) If $w=w_{\bar{r},z}$, then $||w||_{C^2(\bar{\Gw}\cap B_{\frac{\bar r}{2}}(z))}\leq C$ where $C$ is a constant independent of $z$.
\es

\nind{\Remark}: If $g\in\mathcal{G}_0$ and satisfies global KO condition then a g-barrier exists at every point $z\in\bdw$. \\
In fact one can get a stronger implication:
\blemma{global-barrier}
If $\Gw$ is a domain uniformly of class $C^2$ and $g\in\mathcal{G}_0$ satisfies global KO condition in $\Gw$ then $g$ satisfies the global barrier condition.
\es
For proof see \cite[Lemma 3.1.10]{MV}.

\bdef{g-admissible}
Let $\Gw$ be a bounded domain and $g\in\mathcal{G}_0$. A measure $\nu\in\CM(\bdw)$ is called {\it g-admissible}  if
$$g\circ \mathbb{P}_{\Gw}(|\nu|)\in L^1(\Gw,\rho) \ \text{and} \
g\circ \big(-\mathbb{P}_{\Gw}(|\nu|)\big)\in L^1(\Gw,\rho),$$ where
$$\mathbb{P}_{\Gw}(\nu)(x)=\int_{\bdw} P(x,y)d\nu(y),$$ and
$P$ is the Poisson kernel of $-\Gd$ in $\Gw$.
\es

It is known from \cite[Lemma 4.4]{MV1} that if $\Gw$ is a bounded domain, $g\in\CG_0$ and $\nu$ is {\it g-admissible}, then problem \eqref{eq-g} possesses a unique solution.

\vspace{2mm}

It is known from \cite{DG} that, equation \eqref{1'} possesses a global barrier when $g(x,t)=\rho^{\ga}|t|^{q-1}t$ where $\ga>0$ and $q>1$. The notion of global barrier condition (see Definition \ref{d:barrier}) did not appear in their work but the construction of their proof establishes the fact that the nonlinearity $g$ given by $\rho^{\ga}|t|^{q-1}t$, with $\ga>0$ and $q>1$, satisfies global barrier condition. Also note that when $-2<\ga\leq 0$, existence of global barrier follows from the fact that $g$ satisfies global KO condition.
In the next proposition we establish a sufficient condition for the nonlinearity $g\in\mathcal{G}_0$ to satisfy a global barrier condition.

\bdef{C2-char} A possibly unbounded domain $\Gw$ is uniformly of class $C^2$ if it satisfies the following conditions:
\begin{itemize}
\item[(i)]There exists $r_{\Gw}>0$ such that, for every $X\in\bdw$ there exists a set of coordinates $\xi=\xi^X$ and a function $F_X\in C^2(\BBR^{N-1})$ such that $F_X(0)=0$, $\nabla F_X(0)=0$ and
$$\Gw\cap B_{r_{\Gw}}(X)=\{\xi:|\xi|<r_{\Gw},\ \xi_1>F_X(\xi_2,\cdots, \xi_N)\}.$$
\item[(ii)] The set $\{F_X: X\in\bdw\}$ can be chosen so that
$$||\bdw||_{C^2}:=\mbox{sup}\{||F_X||_{C^2(\bar B_{r_{\Gw}}(0))}: X\in\bdw\}<\infty $$
and there exists $\kappa\in C(0,1)$ such that $D^2F_X$ has modulus of continuity $\kappa$ for every $X\in\bdw$. The pair $(r_{\Gw}, ||\bdw||_{C^2})$ is called a $C^2$ characteristic of $\Gw$.
\end{itemize}

\es

\bprop{barrier}
Let $\Gw$ be a domain (not necessarily bounded), uniformly of class $C^2$. We assume that $g\in\mathcal{G}_0$ satisfies {\it local} KO condition and there exists  $C, \ T>0$ such that
\begin{equation}\label{h-infinity}
g(x,t)>C\rho^{\ga_{0}}t^{\Gg}\quad\text{when}\quad t>T,
\end{equation}
 for some $\ga_{0}\geq 0$ and $\Gg>1$.
 Let $z\in\bdw, \ 0<r<\frac{\gb_{0}}{2}$. Then there exists a positive solution $W=W_{z,r}$ of
\eqref{1'} in $B_{r}(z)\cap\Gw$ such that
\begin{equation}\label{12}
\BAL
&(i)\q W\in C(B_{r}(z)\cap\bar\Gw),\quad W=0 \ \text{on} \ \bdw\cap B_{r}(z), \\
&(ii)\q {\lim_{x\to y}}
W(x)=\infty \quad \forall y\in\prt B_{r}(z)\cap\Gw.
\EAL
\end{equation}
Furthermore, if $r=\frac{3\gb_{0}}{4}$, the norm of the solution $W_{z,r}$ in $C^{2}(\bar\Gw\cap B_{\frac{\gb_{0}}{2}}(z))$ is bounded by a constant depending only on $N,\ga_0$ and the $C^2$ characteristics of $\Gw$ but not on $z$.
\es
A proof of the proposition  is provided in the Appendix.

\bdef{maximal-sol}
Let $F$ be a compact subset of $\bdw$ and $g\in\mathcal{G}_0$. Denote by $\mathcal{V}_F$ the family of all non-negative solutions of \eqref{1'} such that
$$u\in C(\bar\Gw\setminus F)\quad\text{and}\quad u=0 \ \text{on} \ \bdw\setminus F.$$
We say that $U$ is the maximal solution relative to $F$ if $U\in\mathcal{V}_F$ and $U$ dominates every solution in $\mathcal{V}_F$.
\es
\blemma{maximal1}
 Suppose that $g\in\mathcal{G}_0$ satisfies the local KO condition. Let $F$ be a compact subset of $\bdw$ such that $\mathcal{V}_F$ is not empty. Then $U_F:=\text{sup}\ \mathcal{V}_F$ belongs to $\mathcal{V}_F$. Thus $U_F$ is the maximal solution relative to $F$.
\es
\begin{proof}
See \cite[Lemma 3.2.3] {MV}.
\end{proof}

\bcor{maximal2}
Let the function $h$ in \eqref{def-H} satisfy \eqref{2} and KO condition.  In addition, assume that \eqref{1} possesses a global barrier and $F$ is a compact subset of $\bdw$ such that $\mathcal{V}_F$ is not empty. Then there exists a positive solution $U_F$ of \eqref{1}
that vanishes on $\bdw\setminus F$ and dominates every solution with this property.
\es
\begin{proof}
Let $u_n$ be the solution of \eqref{1} with $u_n =n \q\text{on}\ F$ and
$u_n=0 \q\text{on}\ \bdw\setminus F$. Since $h$ satisfies KO condition, applying Remark \ref{d:KO} we obtain $u:=\lim_{n\to\infty}u_n$ exists. Let
$z\in\bdw\setminus F$ and $0<r<\frac{\gb_0}{2}$. Then as \eqref{1} possesses a global barrier using Proposition \ref{p:barrier}, we obtain
$$u_n\leq W_{z,r}\q\text{in}\ B_{\frac{\gb_0}{2}}(z)\cap\Gw \q\forall\ n.$$
Thus $u$ vanishes on $\bdw\setminus F$. Now define $\CV_F$ and $U_F$ as in \ref{d:maximal-sol}. Therefore applying Lemma \ref{l:maximal1}, the result follows.
\end{proof}

\section{\bf The similarity transformation and global estimates of the solutions}
A basic tool in our presentation is a \textit{similarity transformation} associated with \eqref{1}, denoted by $T_{a}^{\ga}, \ a>0$,
given by
\begin{equation}\label{similarity}
T_{a}^{\ga}u(x)=a^{\ga}u(ax).
\end{equation}
If $u$ is a weak solution of  \eqref{1} in $\Gw$ then $T_{a}^{\ga}u$ is a weak solution of this equation in $\frac{1}{a}\Gw$. If $\Gw=\frac{1}{a}\Gw$ and $u=T_{a}^{\ga}u$ for every $a>0$, we say that $u$ is a \textit{self-similar solution}.

The next lemma establishes a global estimate of  solutions of \eqref{1}  assuming the Keller-Osserman condition.

\blemma{global}
Let the function $h$ in \eqref{def-H} satisfy the conditions (i) and (ii) in \eqref{2} and KO condition. Then there exists a constant $C=C(N,\ga,h)$ such that for every $\ga>0$ and every solution $u$ of \eqref{1},
\begin{equation}\label{KO-alpha}
|u(x)|\leq C \rho(x)^{-\ga}.
\end{equation}
\es

\begin{proof}
If $u$ is a solution of \eqref{1} then, by Kato's inequality \cite [Proposition 1.5.4] {MV}, $|u|$ is a subsolution; therefore it is sufficient to prove the lemma for $u>0$.

We fix a point $x\in\Gw$ and denote
$$R=\gr(x)/2,\q \Gw^{'}=B_{R}(x).$$
Then
\begin{equation}\label{a1}
  R\leq\rho(y)\leq 3R \forevery y\in\Gw^{'}.
\end{equation}
Using this fact and the monotonicity of $h$ we obtain
\begin{equation}\label{7}
(3R)^{-2-\ga}h(R^{\ga}u)\leq\rho^{-2-\ga}h(\rho^{\ga}u)
\leq R^{-2-\ga}h(3^{\ga}R^{\ga}u) \ \text{in} \ \Gw^{'}.
\end{equation}
Since $u$ satisfies \req{1} and $\Gw'\sbs \Gw$,
\begin{equation}\label{8}
-\Gd u+(3R)^{-2-\ga}h(R^{\ga}u)\leq 0 \quad \text{in} \quad \Gw^{'}.
\end{equation}
Let $v$ be the weak solution of the boundary value problem
\begin{equation}\label{9}
-\Gd v+(3R)^{-2-\ga}h(R^{\ga}v)= 0 \quad \text{in} \quad \Gw^{'},\q v=u\q\text{on }\prt\Gw'.
\end{equation}
By the comparison principle, $u\leq v$ in $\Gw^{'}$.

 Put $c_1:=(3R)^{-2-\ga}$ and $c_2:=R^{\ga}$. Since $h$ satisfies the KO condition, $\tilde h(\cdot):=c_1h(c_2\cdot)$ also satisfies this condition. Denote,
 $$H_1(t)=\int_0^th(s)ds,\q \psi(a)=\int_{a}^{\infty}\frac{ds}{(2 H_1(s))^{\frac{1}{2}}}\forevery a>0,\q \phi=\psi^{-1}.$$
 Define $\tl H_1$, $\tl\psi$ and $\tl\phi$ in the same way with $h$ replaced by $\tl h$.

As $v$ satisfies \req{9},
$$v(y)\leq\tilde\phi(\dist(y,\bdw')\q\text{in }\Gw^{'},$$
(see \cite{K}, \cite{BM}).
 A simple calculation yields  $\tilde\phi(s)=\frac{1}{c_2}\phi(\sqrt{c_1c_2}s)$, $s>0$.  Therefore
\begin{equation}\label{11}
v(x)\leq\tilde\phi(R)\leq\tilde\phi(\frac{\rho(x)}{3})=\frac{1}{c_2}\phi(\sqrt{c_1c_2}\frac{\rho(x)}{3}).
\end{equation}
Here we use \req{a1} and the fact that $\tilde\phi$ is decreasing. Substituting the values of $c_1$ and $c_2$, using again \req{a1} and the fact that $R=\gr(x)/2$ we obtain
$$v(x)\leq C\rho(x)^{-\ga}$$
where  $C=2^\ga\phi(\frac{2}{9}3^{-\frac{\ga}{2}})$. Since $u\leq v$ in $\Gw'$, \req{KO-alpha} follows.
\end{proof}

Assuming existence of barrier at every point of $\bdw$ we can improve the inequality in Lemma \ref{l:global}.

\bprop{global1}
Let $\Gw$ be a domain (not necessarily bounded) uniformly of class $C^2$ and  $F$ be a compact subset of $\bdw$. In addition, assume that the function $h$ in \eqref{def-H} satisfies \eqref{2}  and \eqref{1} possesses a global barrier. Then there exists a constant $C$ depending only on $N,\; \ga,\;h$ and the $C^2$ characteristic of $\Gw$ such that for every solution $u$ of \eqref{1} vanishing on $\bdw\setminus F$,
\begin{equation}\label{29}
|u(x)|\leq C\rho(x)\text{dist}(x,F)^{-\ga-1}
\end{equation}
for every $x\in\Gw$ with dist$(x,F)\leq(1+\gb_0)^{-1}$. If $\Gw$ is bounded then \eqref{29} holds for every $x\in\Gw$.
\es

\begin{proof}
If $u$ is a solution then by Kato's inequality $|u|$ is a subsolution and the existence of barrier guarantees that the smallest solution above $|u|$ vanishes on $\bdw\sms F$. Therefore it is enough to prove the proposition for $u>0$. Let $z\in\bdw\setminus F$ and $\gg(z):=\frac{1}{2}$dist$(z,F)$. Now if $u$ is a solution to \eqref{1} then $u_{z}(x):=\gg(z)^{\ga}u(\gg(z)x)$ is a solution to
\begin{equation*}
-\Gd u_{z}+\rho_\gg(x)^{-2-\ga}h(\rho_\gg(x)^{\ga}u_{z}(x))=0 \quad\text{in}\quad
\frac{1}{\gg(z)}\Gw.
\end{equation*}
where $\rho_\gg(x)=\text{dist}\left(x,\prt\frac{\Gw}{\gg(z)}\right)$

Let $r=\frac{3\gb_{0}}{4}$min$(1,\frac{1}{\gg(z)})$. We assume $\gg(z)<1$ so that $r=\frac{3\gb_{0}}{4}$. Then the solution $W_{z,r}$ mentioned in Proposition \ref{p:barrier} satisfies
\begin{equation}\label{41}
u_z<W_{z,r} \quad \text{in}\quad B_{\frac{\gb_{0}}{2}}(z)\cap\frac{1}{\gg(z)}\Gw.
\end{equation}
Thus $u_z$ is bounded in $B_{\frac{\gb_{0}}{2}}(z)\cap\frac{1}{\gg(z)}\Gw$ by a constant $C$
depending only on $N,\ga$ and the $C^2$ characteristic of $\frac{1}{\gg(z)}\Gw$.
Since $\gg(z)<1$, $C^2$ characteristic of $\Gw$ is also $C^2$ characteristic of $\frac{1}{\gg(z)}\Gw$. Therefore the constant $C$ can be chosen independent of $z$.
Now applying the mean value theorem on $h$ and using Lemma \ref{l:global} we note that, there exists $$c_x\in(0, \rho_\gg(x)^{\ga}u_z)\subset(0,\ C)$$ with
\begin{equation}\label{u_z}
0 =-\Gd u_z+[\rho_\gg(x)^{-2}h^{'}(c_x)]u_z \q\text{in}\ \frac{1}{\gg(z)}\Gw.
\end{equation}
Note that as $h$ is non-decreasing and $C^1$,  $h^{'}(c_x)$ is non-negative and uniformly bounded in $\frac{1}{\gg(z)}\Gw$.
Let $v_z$ denote the solution of
\begin{equation}\label{42}
\BAL
-\Gd v+V(x)v &=0  \quad&&\text{in}
\ \ B_{\frac{\gb_0}{2}}(z)\cap\frac{1}{\gg(z)}\Gw\\
v&=u_z \ \ &&\text{on} \ \prt(\frac{1}{\gg(z)}\Gw\cap B_{\frac{\gb_0}{2}}(z))
\EAL
\end{equation}
where $V(x)=h^{'}(c_x)$ which is non-negative and uniformly bounded. Therefore as $u_z$ is uniformly bounded on $\prt(\frac{1}{\gg(z)}\Gw\cap B_{\frac{\gb_0}{2}}(z))$, applying weak maximum principle we obtain $v_z$ is uniformly bounded in $\frac{1}{\gg(z)}\Gw\cap B_{\frac{\gb_0}{2}}(z)$. Now as $v_z=0$ in $\bdw\cap B_r(z)$, applying Hopf's lemma we obtain
\begin{equation}\label{43}
 v_z(x)\leq C_1 \rho_\gg(x) \ \forall \ x\in B_{\frac{\gb_0}{2}}(z)\cap\frac{1}{\gg(z)}\Gw.
\end{equation}
Also note that, from \eqref{u_z} and \eqref{42} we obtain that $u_z$ is a subsolution of \eqref{42}. Therefore $u_z\leq v_z$ in $B_{\frac{\gb_0}{2}}(z)\cap\frac{1}{\gg(z)}\Gw$. Hence using \eqref{43} and the definition of $u_z$ we obtain
$$u(x)\leq C_1 \text{dist}(x,\bdw){\gg(z)}^{-\ga-1}\ \ \forall \ x\in B_{\gg(z)\frac{\gb_0}{2}}(z)\cap\Gw. $$
Let $x\in\Gw_{\gb_0}$ and assume that
$$\rho(x)=dist(x,\bdw)\leq\gb_0 \text{dist}(x,F), \quad\text{dist}(x,F)<\frac{1}{1+\gb_0}.$$
Let $z$ be the unique point on $\bdw$ such that dist$(x,z)=\rho(x)$. Then
$$2\gg(z)\geq\text{dist}(x,F)-\rho(x)\geq(1-\gb_0)\text{dist}(x,F).$$
Therefore $$u(x)\leq C_1\rho(x)\big(\frac{1}{2}(1-\gb_0)\text{dist}(x,F)\big)^{-\ga-1}.$$
Now if $\rho(x)>\gb_0\text{dist}(x,F)$, by Lemma \ref{l:global}
$$u(x)\leq C\rho(x)^{-\ga}\leq C\rho(x)\big(\gb_0\text{dist}(x,F)\big)^{-\ga-1}.$$
Thus the proposition holds for every $x\in\Gw_{\gb_0}$ such that dist$(x,F)<(1+\gb_0)^{-1}$.
Now if $\Gw$ is bounded, by maximum principle $u$ is bounded in $\{x\in\Gw:\text{dist}(x,F)\geq (1+\gb_0)^{-1}\}$. Therefore by maximum principle \eqref{29} is true in this set and therefore for every $x\in\Gw$.
\end{proof}

\section{\bf Subcriticality}
In this section we consider the boundary value problem \req{1''},
$\nu\in\CM(\bdw)$. $H$ is defined as in \eqref{def-H} and $h$ is assumed to satisfy \eqref{2} and the KO condition. In addition, we assume that \eqref{1} possesses a global barrier and $\Gw$ is a bounded domain.

For $p>1$,  $L^p_w(\Gw,\tau)$ denotes the weak $L^p$ space with the norm
$$||f||_{L^p_w(\Gw,\tau)}=
\sup\displaystyle\left\{\frac{\int_{\gw}|f|d\tau}{\tau(\gw)^\frac{1}{p^{'}}}:\, \gw\subset\Gw, \ \gw\text{ measurable,}\; 0<\tau(\gw)<\infty.\right\}$$
where $\frac{1}{p}+\frac{1}{p^{'}}=1.$

\nind\textit{Notation.} We denote by
$\BBP(\nu)$ the solution of
\begin{equation}\label{28'}
\Gd u =0 \q  \text {in} \q \Gw, \qq u=\nu \q  \text {on} \q \bdw.
\end{equation}

\blemma{continuity}The mapping
$$\BBP:\CM(\bdw)\mapsto L^{\frac{N+\gb}{N-1}}_w(\Gw,\rho^{\gb})\q \forall\gb\geq -1$$
is continuous relative to the norm topologies.
\es
\begin{proof}
See \cite[Lemma 2.3.3(ii)]{MV}
\end{proof}
\bth{uni-int}Let $f\in C(\BBR)$ and be a monotone increasing function such that $$f(t)\to 0 \q\text{as}\q t\to 0$$  $\tau$ be a positive measure in
$\CM(\Gw)$. Assume that, for some $p\in(1,\infty)$, $$\int_{0}^{1}f(r^{-\frac{1}{p}})dr<\infty.$$ If $\{w_n\}$ is a bounded sequence in $L^p_w(\Gw,\tau)$ then $\{f\circ w_n\}$ is uniformly integrable in $L^1(\Gw,\tau)$. Furthermore the modulus of uniform integrability depends only on the bound of the sequence, the $C^2$ characteristic of $\Gw$ and $\diam{\,\Gw}$.
\es
\begin{proof}
See \cite[Thm. 2.3.4]{MV}.
\end{proof}

\noindent {\it Notation.}\hskip 2mm $\CM(\Gw,\rho)$ denotes the space of signed Radon measures
$\mu$ in $\Gw$ such that
$$\rho\mu\in\CM(\Gw)\q\text{where}\q
\rho(x):=\text{dist}\,(x,\prt\Gw).$$
\noindent The norm of a measure $\mu\in\CM(\Gw,\rho)$ is given by
$$||\mu||_{\Gw,\rho}=\int_{\Gw}\rho \ d|\mu|.$$

\vspace{2mm}

\nind{\bf Proof of Theorem \ref{t:continuity}:}  Let $\nu\in\CM(\bdw)$. As $h$ is convex  and odd we have,
\begin{equation}\label{32}
\BAL
\|\rho^{-2-\ga}h\big(\pm\rho^{\ga}\mathbb{P}_{\Gw}(|\nu|)\big)\|_{L^{1}(\Gw, \rho)}&=\int_{\Gw}\rho^{-(1+\ga)}h\big(\rho^{\ga}\int_{\bdw} P(x,y)d|\nu|(y)\big)dx\\
&=\int_{\Gw}\rho^{-(1+\ga)}h\big(\fint_{\bdw}|\nu|(\bdw) \rho^{\ga} P(x,y)d|\nu|(y)\big)dx\\
&\leq\int_{\Gw}\rho^{-(1+\ga)}\displaystyle\left[\fint_{\bdw}h\big(|\nu|(\bdw) \rho^{\ga} P(x,y)\big)d|\nu|(y)\right]dx\\
&=\fint_{\bdw}\int_{\Gw}\rho^{-(1+\ga)}h\big(|\nu|(\bdw) \rho^{\ga} P(x,y)\big)dx d|\nu|(y).
\EAL
\end{equation}

Since $\Gw$ is a bounded domain of class $C^2$,
\begin{equation}\label{P(x,y)}
   c_1^{-1}|x-y|^{1-N}\leq P(x,y)\leq c_1|x-y|^{1-N} \quad\forall (x,y)\in (\Gw\times\bdw),
\end{equation}
where the constant $c_1$ depends only on the $C^2$ characteristic of $\Gw$ and its diameter.

Therefore \req{subcon1} is equivalent to
\begin{equation}\label{subcr1}
\int_\Gw\rho^{-2-\ga}h(c\rho^{\ga}|x-y|^{1-N})\rho\,dx<\infty,
\end{equation}
for every $c>0$ and every $y\in \bdw$. Passing to spherical coordinates
we see that \req{subcr1} implies
$$\int_0^1 s^{-1-\ga}h(cs^{\ga+1-N})s^{N-1}ds<\infty\quad \forall c>0.$$
In fact this inequality is equivalent to \req{subcr1} and therefore to \req{subcon1}. By a substitution of variables this  inequality reduces to \req{subcon2}(a).
Hence, \req{subcr1} implies that, for every $M>0$ there exists a constant $c(M)$, depending only on $M$, the $C^2$ characteristic of $\Gw$ and its diameter, such that
\begin{equation}\label{P-uniform}
\|\rho^{-2-\ga}h(c\rho^{\ga}P(\cdot,y))\|_{L^{1}(\Gw, \rho)}\leq c(M) \quad \forall y\in \bdw,\;\forall c\in [-M,M].
\end{equation}
Finally \req{32} and \req{P-uniform} imply that problem \req{1''} has a solution for every $\nu\in\CM(\bdw)$.

\vspace{2mm}

We turn to the proof of the second assertion of the theorem.
\medskip

\nind {\bf Assertion 1.}\hskip 2mm Let $\{\nu_k\}$ be a bounded sequence in $\CM(\bdw)$ and put $u_k:=\mathbb{P}(\nu_k)$. Then $\{H(\rho, u_k)\}$ is uniformly integrable in $L^{1}(\Gw,\rho)$.
\medskip

We may and shall assume that $\nu_k$ is a positive measure for each $k$.

 Since $\{\nu_k\}$ is bounded in $\CM(\bdw)$ and \eqref{subcon1} implies that $\ga>N-1$ (see \eqref{N-1}),
\begin{equation*}
\BAL
\rho^{\ga}u_k(x) &=\int_{\bdw}P(x,y)\rho^{\ga}(x)d\nu_k(y)\\
&\leq C\int_{\bdw}|x-y|^{-N+1}\rho(x)^{\ga} d\nu_k(y)\\
&\leq C\int_{\bdw}|x-y|^{-N+1+\ga} d\nu_k(y)<C_1.
\EAL
\end{equation*}
Thus $\{\rho^{\ga}u_k\}$ is uniformly bounded in $\Gw$, say by $C_1$.  By assumption
\eqref{h-zero} there exists $M'$ such that
$$h(t)\leq M' t^{1+\ge} \ \forevery\ t\in [0, C].$$ Consequently
$$\gr^{-1-\ga}h(\gr^{\ga}u_k)\leq M'\gr^{-1+\ge\ga}u_k^{1+\ge}.$$
Hence to prove the assertion it's enough to show that $u_k^{1+\epsilon}$ is uniformly integrable in
$L^1(\Gw,\rho^{\epsilon\ga-1})$.  As $\ga>N-1$ we have
$$\frac{(N-1)(1+\ge)}{N+\ge\ga-1}<1.$$
Let $f$ be the function given by $f(s)=s^{1+\ge}$, $s>0$.  Then,
\begin{equation}\label{1.8-1}
   \int_{0}^{1}f(r^{-\frac{N-1}{N+\epsilon\ga-1}})dr<\infty.
\end{equation}

By Lemma \ref{l:continuity},
$\{u_k\}$ is bounded in $L^{\frac{N+\gb}{N-1}}_{w}(\Gw,\rho^{\gb})$ for every $\gb>-1$.
 Choose $\gb=\epsilon\ga-1$. Then, by \req{1.8-1} and
Theorem \ref{t:uni-int} with  $p=\frac{N+\ge\ga-1}{N-1}$ it follows that $\{f(u_k)\}$ is uniformly integrable in
$L^1(\Gw,\rho^{\epsilon\ga-1})$. This proves the assertion.
 \medskip

Now, let $\{\nu_k\}$, $\mu$ and $v_k$ be as in the second part of the theorem.
 As $v_k\leq u_k$ and $\{H(\rho, u_k)\}$ is uniformly integrable (and therefore bounded) in $L^1(\Gw,\rho)$ it follows that $\{H(\rho, v_k)\}$ is bounded in $L^1(\Gw,\rho)$. Hence $\{\Gd v_k\}$ is bounded in this space and \consy $\{v_k\}$ is bounded in $W^{1,p}_{loc}(\Gw)$ for every $p\in[1,\frac{N}{N-1})$. Therefore there exists $v\in L^1(\Gw)$ such that, up to a subsequence, $v_k\to v$ a.e in $\Gw$. This fact and the uniform integrability of $\{H(\rho, v_k)\}$ in $L^1(\Gw,\rho)$ imply that $H(\rho, v_k)\to H(\rho, v)$ in this space. As $\nu_k\rightharpoonup\mu$ it follows that $u_k\to\BBP[\mu]$ in $L^1(\Gw)$. Therefore, since $v_k\leq u_k$ and $v_k\to v$ a.e in $\Gw$, it follows that $v_k\to v$ in $L^1(\Gw)$. These facts together with the weak convergence $\nu_k\rightharpoonup\mu$  imply that $v$ is the weak solution of problem \req{1''} with $\nu=\mu$.

\qed

\noindent\textit{Remark.} The proof of \rth{continuity} actually yields a stronger version of the second part:
\bth{continuity+} Assume that \req{subcon2} and \req{h-zero} hold.
Let $\{\Gw_k\}$ be a sequence of $C^2$ domains with a uniform $C^2$ characteristic \sth
$\Gw_k\sbs B_R(0)$ for some $R>0$ and all $k$. Assume that $\Gw_k\to\Gw$ in the sense that $\dist(\bdw_k,\bdw)\to 0$.

Let $\nu_k\in \CM(\bdw_k)$, $k\in\BBN$. Let  $O$ be a bounded, open neighborhood of $\bdw$ and assume that $\nu_k\rightharpoonup \nu$ weakly relative to $C(\bar O)$, $\nu\in \CM(\bdw)$. Further let $v_k$ denote the solution of
$$\BAL -\Gd v+H(\gr_k,v)&=0 \quad\text{in }\Gw_k,\\
v&=\nu_k  \quad\text{in }\bdw_k,
\EAL$$
where $\gr_k(x)=\dist(x,\bdw_k)$ for $x\in\Gw_k$. Finally let $v$ be the solution of \req{1''}. Considering $v_k$ (resp. $v$) as functions in $B_R(0)$, defined by zero outside  $\Gw_k$ (resp. $\Gw$) we have: (a) $v_k\to v$ in $L^1(B_R(0))$ and uniformly in compact subsets of $\Gw$ and (b) $\{H(\rho, v_k)\}$ is uniformly integrable in the sense that, for every $\ge>0$ there exists $\gd(\ge)>0$ \sth, for every Borel set $E\sbs B_R(0)$,
$$m(E)<\ge\;\Lra\; \int_{E\cap\Gw_k}H(\gr_k,v_k)\,dx<\gd(\ge).$$
\es

\nind{\bf Proof of Theorem \ref{t:removability}:} Since $|u|$ is a subsolution, it's enough to prove the theorem in the case $u\geq 0$. Without loss of generality let us assume that $y=0$.

We claim that
\begin{equation}\label{trace-cond}
H(\rho, u)\in L^1(\Gw,\rho) \  \text{and}\  u\in L^1(\Gw).
\end{equation}
To prove the claim, let $\eta$ be a function in $C^2(\BBR)$ such that
$$0\leq\eta\leq 1,\q \eta(t)=0\ \text{for}\ t<1, \q \eta(t)=1\ \text{for}\ t>2.$$
Further let $\phi$ be the solution of
$$-\Gd\phi=1 \ \text{in}\ \Gw,\q \phi=0 \ \text{on}\ \bdw.$$
Given $\epsilon>0$, set $\gz_\epsilon=\eta(\frac{|x|}{\epsilon})\phi$. Thus $\gz_\epsilon\in C_2(\Gw)$ and vanishes on $\bdw$ and in a neighborhood of origin. Therefore we have
\begin{equation}\label{u-z1}
\int_{\Gw}-u\Gd\gz_\epsilon+H(\rho,u)\gz_\epsilon=0.
\end{equation}
Let $E_\epsilon=\{x\in\Gw:\epsilon<|x|<2\epsilon\}$. Then a straight forward calculation yields that
\begin{equation}\label{u-z2}
\int_{\Gw}u\Gd\gz_\epsilon\leq C_1\int_{E_\epsilon}u\epsilon^{-1}-\int_{\Gw}u(x)\eta(\frac{|x|}{\epsilon})
\end{equation}
where  $C_1$ is a constant independent of $\epsilon$. By Proposition \ref{p:global1}, $u(x)\leq C\rho(x)|x|^{-\ga-1}$. Therefore in $E_\epsilon$, $u(x)\leq C\epsilon^{-\ga}$. Consequently,
\begin{equation}\label{u-z3}
\int_{E_\epsilon}u\epsilon^{-1}\leq C_2\epsilon^{N-\ga-1}\leq C_3
\end{equation}
where $C_3$ is a constant independent of $\epsilon$.
Therefore combining \eqref{u-z1}, \eqref{u-z2} and \eqref{u-z3} we obtain,
$$\int_{\Gw}\big(u(x)+H(\rho,u)\phi\big)\eta(\frac{|x|}{\epsilon})<C_4$$
where $C_4$ is a constant independent of $\epsilon$. Now as $\phi=O(\rho)$, letting $\epsilon\to 0$ to the previous inequality and applying Fatou's lemma we obtain \eqref{trace-cond}. Hence the claim follows.

Now suppose that $u>0$. By Corollary \ref{c:maximal2}, let $U$ be the maximal solution of \eqref{1} vanishing on $\bdw\setminus\{0\}$. Since $U$ satisfies \eqref{trace-cond}, $U$ must have boundary trace $c\gd_0$, for some $c>0$. Clearly $2U$ is a supersolution with the boundary trace $2c\gd_0$. Therefore the largest solution  dominated by $2U$, say $U'$, has the same trace $2c\gd_0$ which is impossible as $U$ is the maximal solution vanishing on $\bdw\setminus\{0\}$.
\hfill{$\square$}

\bcor{sub-al}
Let $\Gw$ be a bounded $C^2$ domain, $H$ be defined as in \eqref{def-H} and $H$ satisfies global barrier condition.  Assume that $h$ satisfies \req{2}, \eqref{h-zero}, KO condition and $\Gd_2$ condition. Then $H$ is subcritical if and only if $\ga>N-1$.
\es
\begin{proof}By Theorem \ref{t:removability}, $\ga\leq N-1$ implies $H$ is not subcritical. On the other hand Theorem \ref{t:continuity} implies that if condition \eqref{subcon1} is satisfied then $H$ is subcritical. Therefore it is enough to prove that  when $\ga>N-1$, condition \eqref{subcon1} always holds.

We will prove this by negation. Suppose there exists a $y\in\bdw$ such that $\int_{\Gw}\rho^{-(1+\ga)}h(\rho^{\ga}P(x,y))dx=\infty$. Since $\rho(x)\leq |x-y|$, by \eqref{P(x,y)} we have
\begin{equation}
\BAL
\infty=\int_{\Gw}\rho^{-(1+\ga)}h(\rho^{\ga}P(x,y))dx &\leq\int_{\Gw}\rho^{-(1+\ga)}
h(c_1\rho^{\ga}|x-y|^{1-N})dx\\
&\leq\int_{\Gw}\rho^{-(1+\ga)}h(c_1\rho^{\ga+1-N})dx\\
&\leq c_2\int_{0}^{a}r^{-1-\ga}h(c_1 r^{\ga+1-N})r^{N-1} dr\\
&\leq c_3\int_{0}^{b}\frac{h(r)}{r^2}dr \ (\text{as}\ \ga>N-1)
\EAL
\end{equation}
where $b=c_1a^{\ga+1-N}$ and $c_3$ is a positive constant. Therefore $\int_{0}^{b}\frac{h(r)}{r^2}dr=\infty$ but this is impossible since we assumed \eqref{h-zero}. Hence condition \eqref{subcon1} holds which implies that $H$ is subcritical if and only $\ga>N-1$.
\end{proof}

\section{\bf Lower estimate of the singular solution with point singularity}

In this section we study \eqref{1} when $H$ satisfies the subcriticality condition \req{subcon1} and $h$  satisfies the KO condition.
\vspace{2mm}

\nind{\bf Proof of Theorem \ref{t:lower-estimate}:} Here we assume that the set of coordinates in $\BBR^N$ is positioned so that $y$ is the origin and the hyperplane $x_1=0$ is tangent to $\bdw$ at $y$ with the positive $x_1$-axis pointing into the domain.

Suppose that $u$ is a positive solution of \eqref{1} and define
$$\hat H(\rho,t):=\frac{H(\rho, t)}{t}\forevery t>0.$$

\nind\textbf{Assertion 1.}\hskip 2mm \textit{Under the assumptions of the theorem, there exists a sequence $\{\xi_n\}\subset\Gw$  converging to the origin such that
\begin{equation}\label{u-s}
u(\xi_n)\rho(\xi_n)^{N-1}\to\infty.
\end{equation}
}

By negation, if the assertion is not valid there exists $R\in (0,1)$ and a constant $C$ \sth
\begin{equation}\label{R-u}
u(x)<C\rho(x)^{1-N} \q \forall x\in B_R(0)\cap\Gw.
\end{equation}
As $h$ is convex, non-decreasing and $h(0)=0$, it follows that $t\mapsto\frac{h(t)}{t}$ is non decreasing on $(0,\infty)$.
Therefore
 $$\hat H(\rho,u)= \frac{\rho^{-2}h(\rho^{\ga}u)}{\rho^{\ga}u}\leq
 \frac{\rho^{-2}h(C\rho^{\ga+1-N})}{C\rho^{\ga+1-N}}
=\frac{1}{C}\rho^{N-\ga-3}h(C\rho^{\ga-N+1})$$
in $Q:=\Gw\cap B_R(0)$.
Thus
$$f_u(t):=\text{sup}_{\Sigma_t\cap B_R(0)}\hat H\leq\frac{1}{C}
t^{N-\ga-3}h(Ct^{\ga-N+1})$$
and by \eqref{sub},
$$\int_{0}^{1}tf_u(t)dt<\infty.$$ Consequently, by \cite[Lemma 3.1.16]{MV}  $u\in L^1(Q)$ and $H(\rho,u)\in L^1(Q,\rho)$. Therefore $0\in\CR(u)$, which contradicts the assumption that $0\in \CS(u)$.

   Next we observe that
\begin{equation}\label{H_u}
\hat H=\frac{\rho^{-2}h(\rho^{\ga}u)}{\rho^{\ga}u}\leq  C_1\rho^{-2} \quad\text{in}\ \Gw
\end{equation}
Indeed by Lemma \ref{l:global}, $\rho^{\ga}u\leq C$ where $C$ depends only on $N,\ga, h$. Therefore \req{H_u} follows from the fact that
$$\sup_{0<t<C}h(t)/t<\infty.$$

Let $B_n=B_{\rho(\gx_n)/2}(\gx_n)$. As $u$ satisfies the equation $-\Gd u+\hat H u=0$, \req{H_u} and the classical Harnack inequality imply that there exists a constant $\bar c$ \sth
\begin{equation}\label{int-Har}
    \sup_{B_n}u\leq \bar c\inf_{B_n}u.
\end{equation}

We assume that $|\gx_n|<\gb_0/4$. Put
$$\BAL\gg_n&=\rho(\gx_n),& \eta_n&=\gs(\gx_n),\\
 D_n&=\{x\in\Gw:\rho(x)>\gg_n\},& V_n&=B_n\cap\Gs_{\gg_n}.\EAL$$
 Then by \req{int-Har} and \req{u-s}
\begin{equation}\label{b_n}
b_n:=\int_{V_n}u\,dS\geq c_1(\Gw)(\sup_{B_n} u)\gg_n^{N-1}\tin.
\end{equation}

Let $f_{n,k}$ be a function on $\Gs_{\gg_n}$ given by,
\begin{equation}\label{fnk}
f_{n,k}=(k/b_n)u\q\text{in }V_n,\q f_{n,k}=0\q\text{in }\Gs_{\gg_n}\sms V_n.
\end{equation}
Then, by \req{b_n},
\begin{equation}\label{fnk-est}
\BAL
&f_{n,k}\leq kc_1^{-1}\gg_n^{1-N}\q \text{and}\\
&f_{n,k}dS_{\,\Gs_{\gg_n}}\rightharpoonup k\gd_0\q\text{weakly relative to } C(\bar\Gw)
\EAL
\end{equation}
where $dS_{\,\Gs_{\gg_n}}$ denotes the surface element on $\Gs_{\gg_n}$.

Let $w_{n,k}$ denote the solution of the \bvp
\begin{equation}\label{wnk}
\BAL -\Gd w+H(\rho,w)&=0 && \text{in } D_n\\
w&=f_{n,k}&& \text{on } \prt D_n.
\EAL
\end{equation}
Given $k>0$ pick $n(k)$ \sth
$$b_n\geq k \q\text{in }V_n \forevery n\geq n(k).$$
Then $f_{n,k}\leq u$ on $\Gs_{\gg_n}$ and \consy
\begin{equation}\label{wnk<u}
w_{n,k}\leq u\q\text{in }D_n \forevery n\geq n(k).
\end{equation}
Further, by \req{fnk-est},
$$P(x,\eta_n)\geq c_2(\Gw)\gg_n^{1-N}\geq c_3(\Gw)f_{n,k}(x)/k \forevery x\in V_n.$$
Therefore, by the maximum principle,
\begin{equation}\label{P>wnk}
(k/c_3)P(x,\eta_n)\geq w_{n,k}(x) \forevery x\in D_n.
\end{equation}

By the argument employed in the proof of \rth{continuity}, the sequences $\{P(\cdot,\eta_n)\} $ and
$\{H(\rho,(k/c_3)P(\cdot,\eta_n))\}$ are uniformly integrable in $L^1(\Gw)$ and $L^1(\Gw;\rho)$  respectively.  In view of \req{P>wnk} this implies that the sequences $\{w_{n,k}\}_{n=1}^\infty$ and
$\{H(\rho,w_{n,k})\}_{n=1}^\infty$ are uniformly integrable in $L^1(\Gw)$ and $L^1(\Gw;\rho)$ respectively. (Here we refer to the extension of $w_{n,k}$ by zero outside $D_n$) By a standard argument, a subsequence of $\{w_{n,k}\}_{n=1}^\infty$ converges locally uniformly in $\Gw$ to a function $w$.
Since $w_{n,k}$ satisfies \req{wnk} and $\{f_{n,k}\}$ converges weakly as stated in \req{fnk-est}, it follows that $w$ is the (unique) solution of the problem
\begin{equation}\label{lim_wnk}
\BAL -\Gd w+H(\rho,w)&=0 && \text{in } \Gw\\
w&=k\gd_0&& \text{on } \bdw,
\EAL\end{equation}
i.e., $w=u_{0,k}$.
Finally, \req{wnk<u} implies that $u\geq u_{0,k}$. As $k$ was arbitrary we obtain
\req{u_infty<u}.
\qed

\section{\bf The very singular solution}
In this section we study \eqref{1} when $H$ satisfies the subcriticality condition \eqref{subcon1}, $h$ satisfies \req{2} and the KO condition and \eqref{1} possesses a global barrier. These conditions will be assumed without further mention.

Let $\mathcal{U}_{y}$ denote the space of positive solutions of \eqref{1} such that
\begin{equation}\label{37}
u\in C(\bar\Gw\setminus\{y\}),\quad u=0 \ \ \text{on} \ \ \bdw\setminus\{y\}.
\end{equation}
Define
\begin{equation}\label{38}
U_{\infty,y}=\text{sup}\{u\in\mathcal{U}_{y}\}.
\end{equation}
 Then by \rcor{maximal2}, $U_{\infty,y}$ is a solution of \eqref{1} and it satisfies \eqref{37}.

\subsection{\bf In half-space}

\bth{half-space}
Let $\Gw=\BBR^N_{+}=\{x_1>0\}$. In addition to the basic conditions, assume that $h$ satisfies \eqref{h-zero}, \eqref{subcon1} and is strictly convex near zero. Then, for $y\in\prt\BBR^N_{+}$, there exists a unique very singular solution $U_{y}$ at $y$,
\begin{equation}\label{59}
U_{y}(x)=r^{-\ga}w(\gs) \quad\text{where}\quad |x-y|=r,\ \ \gs=\frac{x-y}{r}, \ \ \forall \ x\in\BBR^N_{+},
\end{equation}
and $w$ is the (unique) solution to \eqref{52}.
\es
\begin{proof} Without loss of generality we assume that $y=0$.
Let $U_{\infty,0}$ be defined as in \eqref{38} and let $T_{a}^{\ga}, \ a>0$, be the similarity transformation defined in \eqref{similarity}.
Our basic assumptions imply that $U_{\infty,0}$ is a positive solution vanishing on $\bdw\sms \{0\}$. Therefore $T_a^{\ga}U_{\infty,0}$ is again a solution of \eqref{1} in $\BBR^N_{+}$ which vanishes on $\bdw\setminus\{0\}$. Since $u\mapsto T_{a}^{\ga}u$ is an order preserving $1-1$ mapping
from $\CU_0$ onto itself it follows that $T_{a}^{\ga}U_{\infty,0}=U_{\infty,0}$, i.e.,
$U_{\infty,0}$ is self-similar. Therefore
\begin{equation}\label{50}
U_{\infty,0}(x)=a^{-\ga}U_{\infty,0}(\frac{x}{a}) \quad \forall \ x\in\BBR^N_{+},\q \forall a>0.
\end{equation}
This can be rewritten in the form \begin{equation}\label{51}
\begin{cases}
U_{\infty,0}(x)=r^{-\ga}w(\gs) \quad\text{where}\quad |x|=r,\ \ \gs=\frac{x}{r}, \ \ \forall \ x\in\BBR^N_{+}\\
w(\gs)=U_{\infty,0}(\gs) \quad \forall \ \gs\in S^{N-1}_{+}.
\end{cases}
\end{equation}
A straight forward computation yields that $w$ is a solution of \eqref{52}.
Our next aim is to show that \eqref{52} has a unique solution. If $\tilde w$ is another solution then we can easily check that $v=r^{-\ga}\tilde w$ will be a solution of \eqref{1} in $\BBR^N_{+}$. Therefore $v\leq U_{\infty,0}$ and  $\tilde w\leq w$. By a simple computation, \eqref{52} implies
\begin{equation}\label{53}
\BAL
0&=\int_{S^{N-1}_{+}}(\gs\cdot e_1)^{-(2+\ga)}\big[\tilde wh((\gs\cdot e_1)^{\ga}w)-wh((\gs\cdot e_1)^{\ga}\tilde w)\big]\\
&=\int_{S^{N-1}_{+}}(\gs\cdot e_1)^{-2}\big[\frac{h((\gs\cdot e_1)^{\ga}w)}{(\gs\cdot e_1)^{\ga}w}-\frac{h((\gs\cdot e_1)^{\ga}\tilde w)}{(\gs\cdot e_1)^{\ga}\tilde w}\big]w\tilde w.
\EAL
\end{equation}
As $h$ satisfies \eqref{2}, $\frac{h(t)}{t}$ is nondecreasing. Therefore
\eqref{53} implies
\begin{equation}\label{54}
\frac{h((\gs\cdot e_1)^{\ga}w)}{(\gs\cdot e_1)^{\ga}w}=\frac{h((\gs\cdot e_1)^{\ga}\tilde w)}{(\gs\cdot e_1)^{\ga}\tilde w}\quad \forall \ \gs\in S^{N-1}_{+}.
\end{equation}
 As $h$ is strictly convex and $h(0)=0$,  $\frac{h(t)}{t}$ is strictly increasing. Therefore \req{54} implies that $w=\tilde w$.

\vspace{2mm}

Recall that condition \eqref{subcon1} implies that $\ga>N-1$ (see \eqref{N-1}). Therefore, by \req{51},
\begin{equation*}
\BAL
&\|\rho H(\rho, U_{\infty,0})\|_{L^1\big(B_{1}(0)\cap\{x_1>0\}\big)} =
\int_{S^{N-1}_{+}}\int_{0}^{1}{\rho}^{-(1+\ga)}h((\frac{\rho}{r})^\ga\,w(\gs))r^{N-1}dr\,d\sigma\\
&\geq\int_{S^{N-1}_{+}\cap[\rho>r/2]}\int_{0}^{1}{r}^{-(1+\ga)}h((2)^{-\ga}\,w(\gs))r^{N-1}dr\,d\sigma=\infty.
\EAL
\end{equation*}
Thus $U_{\infty,0}$ is a very singular solution at $0$.

By Theorem \ref{t:lower-estimate}, $u_{\infty, 0}$ (see \req{u_infty}) is the smallest very singular solution at $0$. Clearly $u_{\infty, 0}$ is self-similar; therefore  \req{51} holds for $u_{\infty, 0}$ as well. The uniqueness of the solution of \req{52} implies that $u_{\infty, 0}=U_{\infty, 0}$ so that the very singular solution is unique.
\end{proof}

\subsection{\bf In general domain}
\bdef{nt}
Let $\Gw$ be a  $C^2$ domain, $y\in \bdw$ and $u$ a function in $\Gw$. We say that $u(x)\to \ell$ as $x\to y$ non-tangentially if, for every fixed $c>0$,
$$\lim_{\begin{subarray}{c}x\to y\\ -\mathbf{n}_y\cdot (x-y)>c\end{subarray}}u(x)=\ell.$$
\es

\bprop{sphere}
Let $\Gw$ be a  $C^2$ domain, not necessarily bounded and assume that $h$ satisfies all the assumptions of \rth{half-space}.    For $y\in \bdw$ let $u_{y,\infty}$ be defined as in \req{u_infty}.

Then both $U_{y,\infty}$ and $u_{y,\infty}$  satisfy \req{60}, uniformly in compact subsets of $S^{N-1}_+$.
Consequently
\begin{equation}\label{58}
\frac{u_{y,\infty}}{U_{y,\infty}}\to 1 \quad \text{as $x\to y$ non-tangentially in $\Gw$}.
\end{equation}
\es

\proof Without loss of generality we assume that $y=0$ and that the set of coordinates is positioned so that the hyperplane $x_1=0$ is tangent to $\bdw$ at $0$ with the positive $x_1$- axis pointing into the domain.

Let $u$ stand for either $u_{0,\infty}$ or $U_{0,\infty}$.
 For every $a>0$, let $u_{k}^a$ denote the solution of \req{1} in $\Gw^a:=\rec{a}\Gw$ with boundary trace $k\gd_0$, extended to $\RN\sms \{0\}$, by setting it zero outside $\Gw^a\cup\{0\}$. Further denote
$$ u_{\infty}^a=\lim_{k\tin}u_{k}^a$$
and observe that
\begin{equation}\label{limu2}
  u_{\infty}^a= T_a^{\ga}u_{0,\infty},
\end{equation}
$T_a^\ga$ being defined as in \req{similarity}.
Let $\{a_n\}$ be a \seq of positive numbers converging to zero. In view of the local Keller - Osserman condition and the global barrier condition, we can extract a \sseq (still denoted by $\{a_n\}$) \sth
\begin{equation}\label{ukan}
   u_{k}^{a_n}\to v_k,\q  u^{a_n}_{\infty}\to V
\end{equation}
uniformly in compact subsets of $\RN\sms \{0\}$ and $v_k$ and $V$ are solutions of \req{1} in $\RN_+$ that vanish on $\prt \RN_+\sms \{0\}$.

\medskip

\noindent\textit{Assertion 1.}\hskip 2mm $v_k=u_{k}^0$, namely, the solution of \req{1} in $\RN_+$ with boundary trace $k\gd_0$.
\medskip

Let $s>0$ and put $\Gw^{a,s}:=\Gw^a\cap B_s(0)$. Denote by $v^{a,s}_k$ the solution of the problem,
\begin{equation}\label{eq_was}\BAL
   -\Gd v + H(\gr_{a},v)&=0 \quad\text{in }\;\Gw^{a,s},\\
   v&=k\gd_0 \quad\text{on }\;\prt\Gw^{a,s},
\EAL\end{equation}
where
$\gr_{a}(x)=\dist(x,\bdw^{a}).$
Clearly
\begin{equation}\label{ua<uas}
   u^a_{k}\geq v^{a,s}_k \q\text{in }\;\Gw^{a,s}.
\end{equation}
Keeping $s$ fixed, $\Gw^{a,s}\to \RN_+\cap B_s(0)$ as $a\to 0$. By \rth{continuity+},
$$v_{k}^{a_n,s}\to v_k^{0,s}$$
 where $v_k^{0,s}$ is the solution of \req{1} in $\RN_+\cap B_s(0)$ with boundary data $k\gd_0$. Since $v_{k}^{a_n,s}\leq u_k^{a_n}$ it follows that
$v_k^{0,s}\leq v_k$ for every $s>0$. Clearly $v_k^{0,s}\to u^0_k$ as $s\to\infty$. Therefore
\begin{equation}\label{vk>uk}
 u_{k}^0\leq v_k.
\end{equation}
On the other hand, it is easily verified that $v_k\leq u_{k}^0.$
This proves Assertion 1.

Now $V$ is  a solution of \req{1}  in $\RN_+$ which vanishes on $\prt \RN_+\sms \{0\}$. Furthermore,
$V\geq v_k=u_{k,0}^0$  for every $k>0$.  Therefore $V$ is a very singular solution of \req{1} in $\RN_+$. By  \rth{half-space} (uniqueness)  $V=U_0$ . Since the limit is independent of the \seq we conclude that
\begin{equation}\label{uniformconv}
\lim_{a\dar0}u^{a}_{\infty}(x)= U_0(x)=|x|^{-\frac{2+\ga}{q-1}}\gw(x/|x|)
\end{equation}
uniformly in compact subsets of $\RN_+$.

Let $S_1$ be a compact subset of the open half sphere $S^{N-1}_+$ and let $\bar r$ be a positive number \sth
 $$\{x: |x|\leq \bar r, \;x/|x|\in S_1\}\sbs \Gw.$$
 Then  $S_1\sbs \rec{a}\Gw$ for $0<a\leq \bar r$. As the convergence is uniform on $S_1$, \req{uniformconv} implies,
$$\lim_{a\to0}a^{\frac{2+\ga}{q-1}}u(ax)=\gw(x)\q \text{uniformly }\forall x\in S_1$$
which is equivalent to \req{60}.
\qed

The next proposition is a version of the boundary Harnack principle due to Ancona \cite{An-87}. Given $s,s'>0$ denote by $T(s,\gg)$ the cylinder,
\begin{equation}\label{Tsg}
   T(s,\gg)=\{ \xi  =( \xi', \xi  _N) \in   { \mathbb  R}^{N-1}\times { \mathbb  R}:\, |\xi'|  < s,\, -\gg s < \xi  _N< \gg s  \}.
\end{equation}

\bprop{BHP}
Let $s,\gg$ be positive numbers and let $f$ be a Lipschitz function on $\BBR^{N-1}$ with Lipschitz constant $c_f$ such that
$$f(0)=0,\quad  c_f\leq \gg/10.$$
Denote
\begin{equation}\label{D,Gg}\BAL
    D&=\{(\gx',\gx_N)\in T(s,\gg):\,f(\gx')<\gx_N\}\\
    \Gg&=\{(\gx',\gx_N)\in T(s,\gg):\,f(\gx')=\gx_N\}
\EAL\end{equation}

Let $V\in L^\infty_{loc}(D)$ be a non-negative function and $c_V$ a positive constant such that
\begin{equation}\label{Vcond1}
   V(x)\leq c_V\rho_\Gg(x)^{-2}, \q \gr_\Gg(x)=\dist(x,\Gg)\q \forall x\in D
\end{equation}
and denote $L^V:=-\Gd+V$. Let $u$ be a positive $L^V$-harmonic function (i.e. $L^Vu=0$) in $D$ such that $u$ is continuous in $D\cup \Gg$ and $u=0$ on $\Gg$. Denote by $G^V$ the Green kernel for $L^V$ in $D$. Then there exists a constant $c$ depending only on
$N$, $c_V$ and $\gg$ such that
\begin{equation}\label{u/uA}
 c^{-1} s^{N-2}G^V(x,A')\leq  \frac{u(x)}{u(A)}\leq c s^{N-2}G^V(x,A')
 \q \forall x\in D\cap T(s/2,\gg),
\end{equation}
where $A=(0,\cdots,0,\gg s/2)$ and $A'=\frac{3}{4}A$.
\es

\medskip

The following is a consequence of \cite[Theorem 9.1]{An-97}. We use the notation of the previous proposition.

\bprop{2-ge} Let $D$ be as in \rprop{BHP} and let $V$ be a non-negative function in $D$ such that $V\in L^\infty(E)$ for every set $E\sbs D$ that is bounded away from $\Gg$. Further assume that there exist $\ge>0$ and $c(\ge)>0$ such that
\begin{equation}\label{Vcond2}
   V(x)\leq c\rho_\Gg(x)^{-2+\ge},
\end{equation}
Let $G^V$ (resp. $G^0$) denote the Green kernel for $L^V$ (resp. $-\Gd$) in $D$. Then there exists a constant $c_G$ depending only on $N,\gg,\ge$ and $c(\ge)$
such that
\begin{equation}\label{GV-G0}
    c_G^{-1}G^0\leq G^V \leq c_G G^0 \q \text{in }D\cap T(s/2,\gg).
\end{equation}
\es
\medskip
As $G^0(x, A')\thicksim\rho(x)|x-A'|^{-(N-1)}\thicksim s^{-(N-1)}\rho(x)$ in $D\cap T(s/2,\gg)$, combining the previous two propositions we obtain,
\bcor{2-ge} Let $V_1,V_2$ be two non-negative functions in $D$ satisfying the assumptions of \rprop{2-ge}. Let $u_i$ be a positive $L^{V_i}$ harmonic function in $D$  such that $u_i$ is continuous in $D\cup \Gg$ and $u=0$ on $\Gg$, $i=1,2$.  Then there exists a constant $c$, $C$ depending only on
$N$, $\gg$, $\ge$ and $c(\ge)$ such that
\begin{equation}\label{u1/u2}\BAL
&(i) \q && C^{-1}\frac{\rho(x)}{s}\leq\frac{u_i(x)}{u_i(A)}\leq C\frac{\rho(x)}{s} \q\forall x\in D\cap T(s/2,\gg/2),\ i=1,2,\\
&(ii)\q && c^{-1}\frac{u_1(A)}{u_2(A)}\leq  \frac{u_1(x)}{u_2(x)}\leq c \frac{u_1(A)}{u_2(A)}
 \q \forall x\in D\cap T(s/2,\gg).
\EAL
\end{equation}
where $A=(0,\cdots,0,\gg s/2)$.
\es

\nind{\bf Proof of Theorem \ref{t:existence2}:}
We assume that the set of coordinates is positioned so that $y=0$ and the hyperplane $x_N=0$ is tangent to $\bdw$ at $0$ with the positive $x_N-$axis
pointing into the domain.

Let $U_{\infty,0}$ and $u_{0,\infty}$ be as in \rprop{sphere}. By Theorem \ref{t:lower-estimate}, $u_{\infty,0}$ is the minimal very singular solution at $0$. Therefore in order to establish uniqueness of the very singular solution it is enough to show that $U_{\infty,0}=u_{\infty,0}$.

By  Proposition \ref{p:sphere}, for every  $\gb>1$ there exists a constant $c_{\gb}>0$ such that,
\begin{equation}\label{61}
c_{\gb}^{-1}|x|^{-\ga-1}\rho(x)\leq u_{\infty,0}(x)\leq U_{\infty,0}(x)\leq c_{\gb}|x|^{-\ga-1}\rho(x)
\end{equation}
in the truncated cone
$$E_{\gb}:=\{x\in\Gw:|x|\leq\gb\rho(x),\ |x|<r_\Gw\},$$
 where $(r_\Gw,M_\Gw)$ is the $C^2$ characteristic of
$\Gw$. Hence
\begin{equation}\label{62}
u_{\infty,0}\leq U_{\infty,0}\leq c_1(\gb)u_{\infty,0} \quad\text{in}\quad E_{\gb}.
\end{equation}

Put $$\tilde E_\gb=:=\{x\in\Gw:\gb\rho(x)\leq |x|<\tilde s/4\}.$$
Next we show that, for an appropriate choice of $\gb$ and $\tilde s$, inequality  \req{62} (with a different constant) holds in $\tilde E_\gb$ as well. In this part of the proof we make use of \rcor{2-ge}.

Let $f\in C^2(\BBR^{N-1})$ be a function such that $f(0)=0$, $\nabla f(0)=0$ and, for some $s_0\in (0,r_\Gw)$ and $\gg_0\in (0,1)$,
$$\Gw\cap T(s_0,\gg_0)=\{x=(x',x_N):\,|x'|<s_0,\; f(x')<x_N<\gg_0 s_0\}:=D_0.$$
For $s\in (0,s_0]$ denote,
$$\gg(s)=10\sup_{|x'|<s}|\nabla f(x')|.$$
Note that $\gg(s)\downarrow0$ as $s\downarrow0$ and choose $\tl s\in (0,s_0/2)$ such that
 $$\tl\gg:=\gg(\tl s)\leq \gg_0/4.$$

Next we show that, for  $\gb$ sufficiently large
 the following condition holds:
for every $\gz\in \Gg_0:= \bdw\cap B_{\tl s/2}(0)$,
\begin{equation}\label{etaz}\BAL
&(i)\qq &&\eta_\gz:=\gz-\rec{4}\tl\gg |\gz|\mathbf{n}_\gz\in E_\gb,\\
&(ii)\qq &&\tl E_\gb\sbs \cup_{\gz\in \Gg_0} [\gz,\eta_\gz],
\EAL\end{equation}
where $[x,y]$ denotes the linear segment connecting the points $x,y$.

As $\rho(\eta_\gz)= \rec{4}\tl\gg |\gz|$, \req{etaz}(i) is equivalent to
\begin{equation}\label{etaz1}
|\eta_\gz|\leq \gb\tl\gg|\gz|/4 \q\forall\gz\in \Gg_0.
\end{equation}
Since $|\eta_\gz|\leq (1+\tl\gg/4)|\gz|$,\req{etaz1} holds for $\gb\geq 4(1+\tl\gg/4)/\tl\gg$.

 For $x\in \Gw$ such that $\gr(x)\leq \gb_0$, we denote by $\gs(x)$ the nearest point to $x$ on $\bdw$. Condition \req{etaz}(ii) is equivalent to
\begin{equation}\label{etaz2}
\gb\rho(x)\leq |x|<\tl s/4 \;\Lra\; \rho(x)\leq \rec{4}\tl\gg|\gs(x)|.
\end{equation}
Therefore it is sufficient to show that, for $\gb$ sufficiently large,
\begin{equation}\label{etaz''}
\sup_{\tl E_\gb}|x|/|\gs(x)|<\infty.
\end{equation}
If $\{P_k\}$ is a sequence of points in $\tl E_\gb$ such that $|P_k|/|\gs(P_k)|\to \infty$ then $P_k\to 0$. But, if $\gb>1$,
$$\lim_{k\to\infty}\frac{\rho(P_k)^2+|\gs(P_k)|^2}{|P_k|^2}=1$$
and consequently, $\rho(P_k)/|P_k|\to 1$, which is impossible when $\gb>1$.

In continuation we assume that $\gb$ has been chosen so that \req{etaz} holds. For every
$P\in \Gg_1:= \bdw\cap B_{\tl s/4}(0)\sms\{0\}$, denote by $\gx^P$ the Euclidean coordinates centered at $P$, such that the hyperplane $\gx^P_N=0$ is tangent to $\bdw$ at $P$ and the positive $\gx^P_N$ axis is in the direction of $-\mathbf{n}_P$. Let $T^\Gw_P(s,\gg)$ be defined as in \req{Tsg} with $\gx=\gx^P$. Note that, for $a>0$, $P\in \bdw$,
 $$T^{\Gw/a}_{P/a}(s/a,\gg)=\rec{a}T^\Gw_{P}(s,\gg),\q \mathbf{n}_{P/a}=\mathbf{n}_P,$$
where $\mathbf{n}_{P/a}$ denotes the unit normal vector to $\rec{a}\bdw$ at the point $P/a$.
Put
$$ D_P=T^\Gw_P(|P|/2, \tl\gg)\cap\Gw.$$
Let $u$ be a positive solution of \req{1} in $\Gw$ and let $a>0$. Denote $v^a=T_a^\ga u$ where $T_a^\ga$ is the similarity transformation defined in \req{similarity}. Then $v^a$ is a solution of \req{1} in $\rec{a}\Gw$  vanishing continuously on $\rec{a}\bdw\sms\{0\}$. In particular, if $a=|P|/2$ then $w_P:=v^{a}$ satisfies  \req{1} (with $\rho$ replaced by $\rho_{P}(x)=\dist(x, \frac{2}{|P|}\bdw)$) in
$$D'_P=\frac{2}{|P|}D_P=T'_P(1,\tl\gg)\cap\frac{2}{|P|}\Gw$$
where
$$T'_P(1,\tl\gg):=\frac{2}{|P|}T^\Gw_P(|P|/2,\tl\gg)=
T^{\frac{2}{|P|}\Gw}_{\frac{2}{|P|}P}(1,\tl\gg)$$
and $w_P=0$ on $\Gg_P= D'_P\cap\frac{2}{|P|}\bdw$.

Put $V_P=H(\gr_{P},w_P)/w_P$. We claim that $V_P$ satisfies
\begin{equation}\label{Vcond3}
   V_P(x)\leq c\rho_P(x)^{-2+\ge}\q \text{in }D'_P
\end{equation}
with $\ge$ as in \req{h-zero}. Indeed,

$$V_P=\gr_P^{-2}\frac{h(\gr_P^\ga w_P)}{\gr_P^\ga w_P}\leq c\gr_P^{-2}(\gr_Pw_P)^\ge.$$
In view of the uniform barrier condition, $w_P$ is bounded in $D'_P$ by a constant depending only on $N,\ga,h$ and the $C^2$ characteristic of $\Gw$. Therefore the last inequality implies \req{Vcond3}.

Denote

$$w_{1,P}=T_a^\ga u_{0,\infty}, \q a=\frac{|P|}{2}.$$

Note that

$$w_{1,P}(2x/|P|)=(\frac{|P|}{2})^\ga u_{0,\infty}(x)\q \forall x\in \Gw.$$
By \eqref{61} and \req{etaz}(i), $\eta_P=P-\rec{4}\tl \gg|P|\mathbf{n}_P\in E_\gb$ and consequently
$$ w_{1,p}(2\eta_P/|P|)=(\frac{|P|}{2})^\ga u_{0,\infty}(\eta_P)\geq c_1(\gb)(\frac{|P|}{2})^\ga|\eta_P|^{-\ga-1}\rho(\eta_P).$$
 The point $\eta'_P=2\eta_P/|P|$ is in the same position relative to $T'_P(1,\tl\gg)$ as the point $A$ (defined in Proposition 6.4) relative to $T(s,\gg)$. Therefore applying \eqref{u1/u2}(i) we have,
 \begin{equation}\label{w1p}
 w_{1,p}(z)\geq c_P c_1(\gb)(\frac{|P|}{2})^{\ga}|\eta_P|^{-\ga-1}\rho_P(\eta_P) \rho_P(z) \q\forall z\in D'_P\cap T_{\frac{2P}{|P|}}^{\frac{2}{|P|}\Gw}(1/2, \tl\gg/2).
\end{equation}
As $\rho(\eta_P)=\tl\gg|P|$ and $|\eta_P|=\frac{|P|}{2}|\eta'_P|\thicksim\frac{|P|}{2} |z|$ in $D'_P\cap T_{\frac{2P}{|P|}}^{\frac{2}{|P|}\Gw}(1/2, \tl\gg/2)$, we obtain from \eqref{w1p} that $$w_{1,p}(z)\geq c_P c_2(\gb)|P|^{\ga}\big(\frac{|P|}{2}|z|\big)^{-\ga-1}|P|\tl\gg\rho_P(z)\q\forall z\in D'_P\cap T_{\frac{2P}{|P|}}^{\frac{2}{|P|}\Gw}(1/2, \tl\gg/2).$$
i.e. $$w_{1,p}(z)\geq c_3(\gb)|z|^{-\ga-1}\rho(z)\q\forall z\in D'_P\cap T_{\frac{2P}{|P|}}^{\frac{2}{|P|}\Gw}(1/2, \tl\gg/2).$$
Since $\Gw$ is uniformly of class $C^2$, we may choose $c_P$ to be independent of $P\in\Gg_1$. Hence there exists a constant $c_4(\gb)$ such that,
$$u_{0,\infty}(x)\geq c_4(\gb)|x|^{-\ga-1}\rho(x) \q\forall x\in D_P\cap T_P(|P|/4, \tl\gg/2).$$
This inequality and \eqref{etaz}(ii) imply that
$$u_{0,\infty}(x)\geq c_5(\gb)\rho(x)|x|^{-\ga-1} \q\forall x\in \frac{1}{2}\tl E_{\gb}.$$
Combining this inequality with \eqref{61} we conclude that there exists $r$,$c$ positive such that
$$u_{0,\infty}(x)\geq c\rho(x)|x|^{-\ga-1} \q\forall x\in B_{r}(0)\cap\Gw.$$
Finally combining this inequality with Proposition \ref{p:global1} we obtain
\begin{equation}\label{up-low-esti}
c\rho(x)|x|^{-\ga-1}\leq u_{0,\infty}(x)\leq U_{0,\infty}(x)\leq C\rho(x)|x|^{-\ga-1}
\end{equation}
for every $x\in B_{r}(0)\cap\Gw$.\\

\vspace{2mm}

Claim: \eqref{up-low-esti} holds for every $x\in\Gw$.

\vspace{2mm}

Indeed if $x\in\bar\Gw\setminus B_r(0)$, then $|x|>r$. Now choose $z\in\bdw\setminus B_r(0)$ and consider the cylinder $T(\frac{r}{2}, \tl\gg)$ centered at $z$. Then the point $\eta_z=z-\frac{r\tl\gg}{2}{\bf n}_z$ is in the same position relative to $T(\frac{r}{2}, \tl\gg)$ as the point $A$(defined in Proposition \ref{p:BHP}) relative to $T(\frac{s}{2},\tl\gg)$. Therefore by Corollary \ref{c:2-ge} (in particular \eqref{u1/u2}(i) applied in $D\cap T(\frac{r}{2}, \frac{\tl\gg}{2})$) we obtain,
\begin{equation}\label{hopf type}
u_{0,\infty}(x)\geq c_z \frac{\rho(x)}{r}u_{0,\infty}(\eta_z) \q\forall x\in D\cap T(\frac{r}{2}, \frac{\tl\gg}{2}).
\end{equation}
As $\Gw$ is a bounded $C^2$ domain, we may choose $c_z$ independent of $z\in\bdw\setminus B_r(0)$. Also note that for every $z\in\bdw\setminus B_r(0)$, corresponding $\eta_z$ belongs to $\Gg_{\frac{r\tl\gg}{2}}:=\{x\in\Gw:\rho(x)=\frac{r\tl\gg}{2}\}$. Therefore by maximum principle, $$\inf_{z\in\bdw\setminus B_r(0)}u_{0,\infty}(\eta_z)\geq c_1(r)>0.$$ As $\Gw$ is bounded, we can cover $\{x\in\Gw\setminus B_r(0):\rho(x)\leq\frac{r\tl\gg}{2}\}$  by a finite number of cylinders. Hence from \eqref{hopf type} we conclude that there exists a constant $c=c(r)>0$ such that
$$u_{0,\infty}(x)\geq c(r)\rho(x) \q\forall x\in\Gw\setminus B_r(0), \q \rho(x)\leq\frac{r\tl\gg}{2}.$$ By above inequality and Proposition \ref{p:global1},
\begin{equation}\label{hopf-type2}
c(r)\rho(x)\leq u_{0,\infty}(x)\leq U_{0,\infty}(x)\leq C(r)\rho(x) \q\forall x\in\Gw\setminus B_r(0).
\end{equation}

Therefore \eqref{up-low-esti} with a constant $C$ depending on the parameters mentioned in assertion (ii) of Theorem \ref{t:existence2} holds for every $x\in\Gw$. Hence there exists a positive constant $c$ such that
$$u_{0,\infty}\leq U_{0,\infty}\leq c u_{0,\infty} \q \text{in}\ \Gw.$$

Therefore, as $h$ is convex, a standard argument introduced in \cite{MVsub}, implies that
$$ U_{0,\infty}(x)=u_{0,\infty}(x).$$

\qed

\section{\bf The generalized boundary value problem}

In this section we study the generalized boundary value problem:
\begin{equation}\label{72}
\BAL
-\Gd u+H(\rho, u) &=0,  \quad\text{in} \ \Gw, \\
u &=\nu  \quad\text{on} \ \bdw,\\
u&>0\quad\text{in} \ \Gw,
\EAL
\end{equation}
where $\nu\in\mathcal{B}_{reg}(\bdw)$, $H(\rho,t)$ is given by \eqref{def-H} and satisfies \eqref{subcon1} and $\Gw$ is a bounded $C^2$ domain. Our goal is to prove the existence and uniqueness of the solution for this problem as stated in Theorem \ref{t:existence3}.

\vspace{2mm}

\nind{\bf Proof of Theorem \ref{t:existence3}:}
Existence follows from  \cite[Theorem 3.3.1]{MV} (also see \cite[Theorem 4.16]{MVsub2}) because, in the subcritical case, conditions
(i) and (ii) are satisfied by any measure in $\CB_{reg}(\bdw)$.

Define $F:=S_{\nu}$. By \cite[Theorem 3.3.1]{MV} it is enough to prove the uniqueness result in the case $\nu=\infty\CX_{F}$.

Next we show that, for every compact set $F\subset\bdw$, the maximal solution $U_F$ is the
unique solution of problem \ref{72} with trace $\nu=\infty\CX_{F}$ . We begin by constructing the minimal solution with this trace.

 Let $\{x_n\}\subset F$ be a sequence dense in $F$ and put $$\nu_k=k\sum_{i=1}^{k}\gd_{x_i}.$$
Let $u_k$ be the unique solution of \eqref{1} with boundary trace $\nu_k$. Thus the sequence $\{u_k\}$ is increasing and since $h$ satisfies  KO condition, using Remark \ref{d:KO} we obtain, $\{u_k\}$ is uniformly bounded in compact subset of $\Gw$. Therefore $$V_F:=\lim_{k\to\infty}u_k$$ is a solution of \eqref{1} which is $\infty$ on the set $ \{x_n\}_{i=1}^n$ and vanishes on $\bdw\setminus F$. Now if $U_{x_n}$ is the unique very singular solution at $x_n$, then
$$V_F\geq U_{x_n}, \quad n=1,2,\ldots.$$ Therefore $x_n\subset S(V_F)\subset F$, now as $S(V_F)$ is closed it follows that $S(V_F)=F$. Therefore $V_F$ is a solution to the problem \eqref{72} with boundary trace $\nu=\infty\mathcal{X}_F$.

Now if $u$ is any positive solution of \eqref{1}
such that $S(u)=F$, then by Theorem \ref{t:lower-estimate} we obtain $$u\geq U_{x}, \quad \forall \ x\in F.$$ Hence $u\geq u_k$ and that implies $u\geq V_F$. Therefore $V_F$ is the minimal solution with boundary trace $\infty\mathcal{X}_F$. It remains to prove that if $U_F$ is the maximal solution vanishing on $\bdw\setminus F$
then $$U_F=V_F.$$
Let $x\in\Gw$ and dist$(x,F)\leq\frac{\gb_0}{2}$, where $\gb_0$ is defined as in the beginning of Section 2. Now choose $y\in F$ such that $|x-y|=$dist$(x,F)$.
Now as $U_y$ is the unique very singular solution at $y$, we have $U_y\leq V_F$. Therefore by \eqref{asymp-esti} we obtain,
$$C^{-1}|x-y|^{-\ga-1}\rho(x)\leq U_y\leq V_F.$$
By Proposition \ref{p:global1} we also obtain,
$$V_F(x)\leq U_F(x)\leq C_1\text{dist}(x,F)^{-\ga-1}\rho(x).$$
Hence there exists a positive constant $C_2$ (depending only on $N,\ga, h$, $C^2$ characteristic of $\Gw$) such that
\begin{equation}\label{75}
C_2^{-1}\text{dist}(x,F)^{-\ga-1}\rho(x)\leq V_F(x)\leq U_F(x)\leq C_2 \text{dist}(x,F)^{-\ga-1}\rho(x)
\end{equation}
for every $x\in\Gw$ such that dist$(x,F)<\frac{\gb_0}{2}$.
 By the same  argument that was used to prove \eqref{hopf-type2}, we deduce that \req{75}, with possibly a larger constant, is valid in the entire domain.
Thus there exists a constant $c>1$ such that
$$V_F\leq U_F\leq c V_F \quad\text{in}\ \Gw.$$
As before this implies $U_F=V_F$.
\qed

\section{\bf Appendix}

\nind{\bf Proof of Proposition \ref{p:barrier}:}
Define $$\tilde g(x,t):=\text{max}\{g(x,t), C\rho^{\ga_0}t^\Gg\}.$$ Thus
\begin{equation}\label{g-comp1}
\tilde g(x,t)\geq C\rho^{\ga_0}t^{\Gg} \quad \forall \ t.
\end{equation}
From the result of Du and Guo \cite{DG}, it is known that \eqref{1'} possesses a global barrier
$W$,  which satisfies (i) and (ii)  in Proposition \ref{p:barrier},
when $g(x,t)=\rho^{\ga}t^q$, $\ga>0$ and $q>1$. Therefore, thanks to \eqref{g-comp1}, it is easy to check that $-\Gd u+\tilde g\circ u=0 \quad\text{in} \ \Gw$, possesses a global barrier. Let $z\in\bdw$, $0<r<\frac{\gb_0}{2}$ and $\tilde V_M$ be the solution of
\begin{equation}\label{A1}
-\Gd u+\tilde g\circ u=0 \quad\text{in} \ \Gw\cap B_r(z),
\end{equation}
with
\begin{equation}\label{A2}
\tilde V_M=
\begin{cases}
0 \quad &\text{on} \quad \bdw\cap B_r(z)\\
M \quad &\text{on}\quad \Gw\cap\prt B_r(z).
\end{cases}
\end{equation}
In addition, assume that $V_M$ be the solution of
\begin{equation}\label{A3}
-\Gd u+ g\circ u=0 \quad\text{in} \ \Gw\cap B_r(z),
\end{equation}
where $V_M$ satisfies the boundary data \eqref{A2}. Now note that
$$g(x,t)\leq\tilde g(x,t)\leq g(x,t)+g_1(x,t),$$
where $g_1(x,t)=C\rho^{\ga_0}\text{min}\{t,T\}^{\Gg}$.
Hence $W_M\leq\tilde V_M\leq V_M$, where $W_M$ is the solution of
\begin{equation}\label{A4}
-\Gd u+(g+g_1)\circ u=0 \quad\text{in} \ \cap B_r(z)
\end{equation}
and $W_M$ satisfies \eqref{A2}. Let $W_0^M$ be the solution of
\begin{equation}\label{A5}
-\Gd W_0^M={g_1}\circ W_M \quad\text{in} \ \Gw\cap B_r(z); \quad  W_0^M=0 \quad \text{on} \  \prt\big(\Gw\cap B_r(z)\big).
\end{equation}
Note that, $g_1$ is bounded in $\Gw\cap B_r(z)$, which implies $g_1\circ W_M$ is uniformly bounded. Hence there exists $C_1>0$, independent of $M$, such that $W_0^M<C_1$.
Combining \eqref{A4} and \eqref{A5} we obtain
\begin{equation}\label{A6}
-\Gd(W_M+W_0^M)+ g\circ(W_M+W_0^M)\geq 0 \quad\text{in} \ \Gw\cap B_r(z),
\end{equation}
and
\begin{equation*}
\tilde W_M+W_0^M=
\begin{cases}
0 \quad &\text{on} \quad \bdw\cap B_r(z)\\
M \quad &\text{on}\quad \Gw\cap\prt B_r(z).
\end{cases}
\end{equation*}
Hence $V_M\leq W_M+W_0^M\leq W_M+C_1$. As $g$ satisfies local KO condition, $\tilde g$ satisfies the same. Therefore passing to the limit we obtain, $\lim_{M\to\infty}\tilde V_M$ exists and the limit is a barrier for the equation $-\Gd u+\tilde g\circ u=0 \quad\text{in} \ \Gw$ at $z$. Now as $W_M\leq\tilde V_M$, we obtain, $\lim_{M\to\infty}W_M=W$. Finally, as $V_M\leq W_M+C_1$, we can conclude that $\lim_{M\to\infty}V_M=V$ is a barrier of \eqref{1'} at $z$.
\hfill{$\square$}

\end{document}